\documentclass{amsart}

\usepackage{float}
\usepackage{amsfonts}
\usepackage{graphicx}
\usepackage{amssymb}
\usepackage{microtype}
\usepackage{tikz}
\usetikzlibrary{arrows}
\usepackage{booktabs}
\usepackage{epstopdf}
\def\MR#1{\href{http://www.ams.org/mathscinet-getitem?mr=#1}{MR#1}}

\theoremstyle{definition}

\numberwithin{figure}{section}
\numberwithin{equation}{section}
\numberwithin{table}{section}

\newcommand{\n}{\eta}

\newcommand{\cal}{\mathcal}

\title{Hidden Multiscale Order in the Primes}

\author{Salvatore Torquato}

\address{Department of Chemistry, Department of Physics,
Princeton Institute for the Science and Technology of
Materials, and Program in Applied and Computational Mathematics, Princeton University,
Princeton NJ 08544}

\author{Ge Zhang}

\address{Department of Chemistry, Princeton University,
Princeton NJ 08544 }

\author{Matthew De Courcy-Ireland}

\address{Department of Mathematics, Princeton University,
Princeton NJ 08544}

\begin{document}

\bibliographystyle{amsplain} 

\begin{abstract}
We study the {pair correlations between} prime numbers in an interval $M \leq p \leq M + L$ with $M \rightarrow \infty$, $L/M \rightarrow \beta > 0$. By analyzing the \emph{structure factor}, we prove, conditionally on the {Hardy-Littlewood conjecture on prime pairs}, that the primes are characterized by unanticipated multiscale order. Specifically, their limiting structure factor 
is that of a union of an infinite number of periodic systems and is characterized by dense 
set of Dirac delta functions.  Primes in dyadic intervals are the first examples of what we call  {\it effectively limit-periodic} point configurations.
This behavior implies anomalously suppressed density fluctuations compared to uncorrelated (Poisson) systems at large length scales, which is now known as hyperuniformity. Using a scalar order metric $\tau$ calculated from the structure factor, we identify a transition between the order exhibited when $L$ is comparable to $M$ and the uncorrelated behavior when $L$ is only logarithmic in $M$. Our analysis for the structure factor leads to an algorithm to  reconstruct primes in a dyadic interval with high accuracy.

\end{abstract}

\maketitle

\section{Introduction} \label{sec:intro}

While prime numbers are deterministic, by some probabilistic descriptors, they can be regarded as pseudo-random in nature.
Indeed, the primes can be difficult to distinguish from a random configuration of the same density. For example, assuming a plausible conjecture, Gallagher proved that the gaps between primes follow a Poisson distribution \cite{gallagher1976distribution}. Thus the number of $M \leq X$ such that there are exactly $N$ primes in the interval $M < p < M + L$ of length $L \sim \lambda \ln{X}$ is given asymptotically by
\begin{equation*}
\# \left\{M \leq X \ \text{such \ that} \  \pi(M+L)-\pi(M) = N  \right\} \sim X \frac{e^{-\lambda}\lambda^{N}}{N!}.
\end{equation*}
Note that Gallagher's interest was in short intervals whereas we analyze those in which the length $L$ is comparable to the lower endpoint $M$. Our primary observation is that for these longer intervals, the primes are highly correlated and ordered on multiple length scales and hence are drastically different from a Poisson distribution. This is demonstrated by the identification of sharp peaks in the \emph{structure factor} of the primes and by large values of the order parameter $\tau$, both of which we define in Section ~\ref{sec:definitions}. In particular, we use the structure factor to detect a large-scale order known as hyperuniformity \cite{To03a}, very different from the uncorrelated behavior one sees in short intervals.

To study different ranges of primes, it is important to take account of the fact that primes become increasingly sparse in longer intervals.
Let $\pi(x)$ denote the {\it prime counting function}, which gives the number of primes
less than $x$. According to the prime number theorem \cite{hadamard1896distribution}, the prime counting
function in the large-$x$ asymptotic limit is given by
\begin{equation}
\pi(x) \sim \frac{x}{\ln(x)} \qquad (x \rightarrow \infty).
\end{equation}
The prime number theorem means that for sufficiently large $x$, the probability that a randomly
selected integer not greater than $x$ is prime is very close to $1 / \ln(x)$,
which can be viewed as a position-dependent number density $\rho(x)$ (number of primes up to $x$ divided
by the interval $x$). This implies that  the primes become sparser as $x$ increases
and hence constitute a {\it statistically inhomogeneous} set of points that are located
on a subset of the odd integers. This simple observation requires that one carefully choose the interval over which the primes are sampled and characterized in order to obtain meaningful results that in general will depend on the chosen interval. If $L$ is much larger than $M$, the density $1/\ln(n)$ drops off appreciably as $n$ ranges from $M$ to $M+L$, and then the system is the very opposite of hyperuniform. On the other hand, $\ln(M+L) = \ln(M) + \ln(1 + L/M)$ is asymptotic to $\ln(M)$ as long as $L/M$ is bounded above. In this case, one can treat the primes as homogeneous with constant density {$1/\ln(M)$. For this paper, we take $L \sim \beta M$ of the same order as $M$ or sometimes smaller to compare with Gallagher's regime.\\

The plausible conjecture Gallagher assumed is a version of the Hardy-Littlewood $m$-tuples conjecture (Theorem X1, p. 61 of \cite{Ha23}). If $\mathcal{H} = (h_1,\ldots ,h_m)$ is a $m$-tuple of integers, then the conjecture gives the number of $n \leq X$ such that all of the shifts $n+h_1,\ldots, n+h_m$ are prime as
\begin{equation} \label{eq:tuples}
\# \left\{n \leq X \ \text{such \ that} \ n+h_j \ \text{all \ prime} \right\} \sim \mathfrak{S}(\mathcal{H})\frac{X}{(\ln{X})^m},
\end{equation}
where
\begin{equation}
\mathfrak{S}(\mathcal{H}) = \prod_{p} \left(1 - \frac{1}{p}\right)^{-m}\left(1 - \frac{\nu_{\mathcal{H}}(p)}{p}\right) 
\label{HLConstant}
\end{equation}
\begin{equation}
\nu_{\mathcal{H}}(p) = \# \left\{\text{distinct} \ h_j \ \text{mod} \ p \right\}.
\end{equation}
When $m=1$, say $\mathcal{H} = \{ h_1 \}$, (\ref{eq:tuples}) simply counts primes less than $X$ (or, strictly, less than $X - h_1$). In (\ref{HLConstant}), since every $\nu(p)$ is 1, one has $\mathfrak{S}(\mathcal{H}) = 1$. Thus the case $m=1$ is the prime number theorem, and it is the only one so far to be proved. Gallagher used all values of $m$ in order to compare the empirical moments of primes in short intervals with the moments of the Poisson distribution. 

Our study of the structure factor ultimately leads to an equivalent formulation of the case $m=2$.
The Hardy-Littlewood constant $\mathfrak{S}(\mathcal{H})$ can be understood as a correction 
to the prediction one would make by imagining that all of the shifts $n+h_{j}$ are prime 
independently with probability $1/\ln(X)$; see \cite{cherwell} and, for this and other senses in which 
the random model fails, \cite{pintz2007}. To summarize the interpretation, note that, for each $p$, 
$(1-1/p)^m$ is the naive chance that each of the shifts would be indivisible by $p$. However, these 
constraints are not independent, and $\nu_{\mathcal{H}}(p)$ is exactly the number of residue classes 
modulo $p$ which $n$ must avoid or else one of the numbers $n+h_j$ would have $p$ as a factor. 
Thus $\mathfrak{S}(\mathcal{H})$ cancels the incorrect guess $(1-1/p)^{m}$ and replaces it with the 
correct $(p - \nu(\mathcal{H}))/p$. The argument advanced by Hardy-Littlewood, which we outline in Appendix A, is however of an altogether different nature.

Probabilistic methods to treat the primes have yielded fruitful insights
about them \cite{Gr95}. Furthermore, there are computationally quick {\it stochastic}
ways to find large primes \cite{miller_1976,rabin_1980, pomerance1980pseudoprimes, baillie1980lucas, atkin1993elliptic}. 
On the other hand, it is known that primes contain unusual patterns, and hence their distribution is not purely random.
Chebyshev observed (circa $1853$) that primes congruent to $3$ modulo $4$ seem to predominate 
over those congruent to $1$. Assuming a generalized Riemann hypothesis, Rubinstein and Sarnak ~\cite{rubinstein_1994}
exactly characterized this phenomenon and more general related results. A computational study on the 
Goldbach conjecture demonstrates a  connection based on a modulo $3$ geometry between the set of even integers 
and the set of primes~\cite{martelli_2013}. 
In $1934$, Vinogradov proved that every sufficiently large odd integer is the sum of three primes~\cite{vinogradov_1937}. 
This method has been extended to cover many other types of patterns~\cite{green_2008,green_2006,green_2008_2,tao_2011}. 
Recently it has been shown that there are infinitely many pairs of primes with some finite gap~\cite{zhang_2014}
and that primes with decimal expansion ending in $1$ are less likely to be followed by another prime ending in $1$~\cite{lemke_2016}.
There is numerical evidence for patterns in the distribution of gaps between primes when these are divided into 
congruence families ~\cite{maynard_2015,dahmen_2001,wolf_1996}.

The present paper is motivated by {certain remarkable} properties of the 
Riemann zeta function $\zeta(s)$, which is a function of a complex variable $s$
that is intimately related to the primes. The zeta function has many different representations,
one of which is the well-known series formula
\begin{equation}
\zeta(s)=\sum_{n=1}^{\infty} \frac{1}{n^s},
\end{equation}
which converges for $Re(s)>1$.  However, $\zeta(s)$ has a unique analytic continuation to the entire complex plane,
excluding the simple pole at $s=1$.
According to the {\it Riemann hypothesis}, the nontrivial
zeros of the zeta function lie along the {\it critical line} $s=1/2 + it$ with $t\in\mathbb{R}$ in the complex plane
{and hence form a one-dimensional point process.}
The nontrivial zeros tend to get denser the higher on the critical line. When the spacings of the zeros are appropriately
normalized so that they can be treated as a homogeneous point process at unity density, the resulting
pair correlation function $g_2(r)$ takes on the simple form $1-\sin^2(\pi r)/(\pi r)^2$ \cite{Mon73}.
This has consequences for the distribution of the primes in short intervals \cite{Mo04}. The corresponding
structure factor $S(k)$ [essentially the Fourier transform of $g_2(r)$] s given by
\begin{eqnarray}
S(k) = \left\{
\begin{array}{lr}
\frac{\displaystyle k}{\displaystyle 2\pi}, \quad 0 \le k \le 2\pi\\\\
\displaystyle{1}, \quad k > 2\pi.
\end{array}\right.
\label{spec}
\end{eqnarray}
Remarkably, this exactly matches the structure factor of the eigenvalues of a random matrix in the Gaussian unitary ensemble
\cite{Dy62a,Me91,Rud96}; see Fig. \ref{nontrivial}.
We see that structure factor goes to zero linearly in $k$ as wavenumber goes to zero and is equal to unity
for $k > 2\pi$. This implies that the normalized Riemann zeros possess an unusual type of correlated
disorder at large length scales known as hyperuniformity~\cite{To03a,To08c}. A hyperuniform point configuration is one in which
$S(k)$ tends  to zero as the wavenumber $k$ tends to zero  ~\cite{To03a}. In such systems, density fluctuations are
anomalously suppressed at very large length scales, a ``hidden"
order that imposes strong global structural constraints.
All structurally perfect crystals and quasicrystals are hyperuniform,
but typical disordered many-particle systems, including gases, liquids, and glasses, are not. Disordered
hyperuniform many-particle systems are exotic states of
amorphous matter that have attracted considerable recent attention \cite{Do05d,To08c,Ba08,Za09,Fl09b,Za11a, Ji14,Ja15,To15,Ma15,Ma16,Go17,He17b,To18a}.

\begin{figure}[bthp]
\includegraphics[  width=3in,keepaspectratio,clip=]{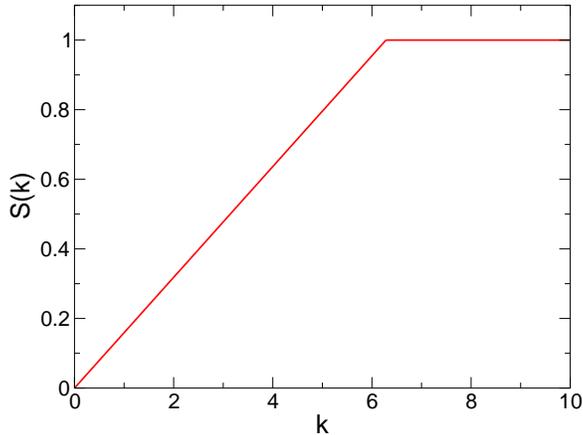}
\caption{Structure factor of the normalized nontrivial zeros of the Riemann zeta function as a function
of the wavenumber [{\it cf.} (\ref{spec})].  This is a special case of a hyperuniform point configuration~\cite{To03a,To08c}.}
\label{nontrivial}
\end{figure}

Because information about the  primes can in principle be deduced  from information about 
the nontrivial  zeros of the zeta function via explicit formulas  \cite{Da80,Te95,Iw04},
one might expect the primes to encode hyperuniform correlations seen in the Riemann zeros. For example,
von Mangoldt's explicit formula for a weighted counting function $\psi(x)= \sum_{p^n < x} \ln(p)$ is given by
\begin{equation}
\psi(x) = x - \sum_{s} \frac{x^{1/2+i\gamma}}{\frac{1}{2} + i\gamma} -\frac{1}{2} \ln(1-x^{-2}) -\ln(2\pi),
\end{equation}
for $x >1$ and $x$ not a prime or prime power (where $\psi(x)$ would have a jump discontinuity). Here, $1/2 + i\gamma$ denotes a nontrivial zero of $\zeta(s)$, meaning that it lies in the critical strip $0 < \text{Re}(s) < 1$. Assuming the Riemann Hypothesis, $\gamma$ is real and the zeros thus form a one-dimensional point process. In any case, allowing for complex $\gamma$, the explicit formula applies unconditionally. The trivial zeros $-2, -4, -6, \ldots$ contribute $\ln(1-x^{-2})$. The explicit formula may be thought of as a Mellin transform of the prime numbers. The structure factor $S(k)$, which is the basis of our investigation, is also a Fourier-type transform of the primes (without any weights) but in a more direct sense:
 \begin{equation}
                   S(k) = \frac{1}{N} \left| \sum_p e^{ikp} \right|^2
                   \label{Skdef1}
                  \end{equation}
where $p$ runs over the primes in the interval $[M,M+L]$, the number of which we denote by $N$. A weighted version of the inner sum, which weights each prime $p$ and also its higher powers $p^{l}$ by $\ln{p}$, has been much studied in connection with the circle method (see section 25 of \cite{Da80}, for example). The behavior of $S(k)$ for small values of $k$ reflects the large-scale correlations between primes.

\begin{figure}[h]
\includegraphics[width=0.85\textwidth,clip=]{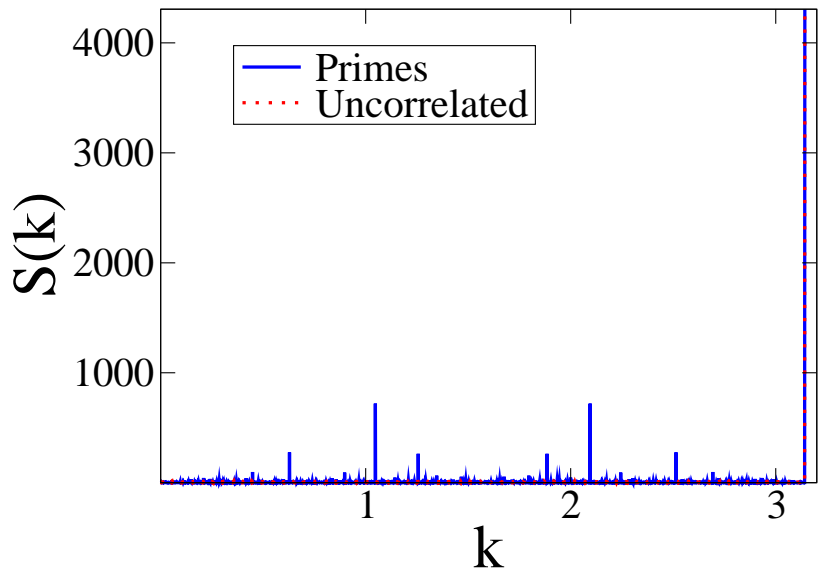}
\caption{The structure factor $S(k)$ of the  prime numbers for $M=10^{10}+1$ and $L=10^{5}$ obtained
in a separate numerical study \cite{Zh18}. It is seen that
it  contains many well-defined Dirac-delta-function like (Bragg-like) peaks of various intensities
characterized by a type of self-similarity. Included in the figure is the corresponding structure factor
for the uncorrelated lattice (Poisson) gas on the integer lattice, whose intensities are barely perceptible on the scale of this figure.}
\label{Sofk-numerical}
\end{figure}

In a very recent numerical study \cite{Zh18}, we and Martelli  examined the pair statistics of the primes, especially the structure
factor $S(k)$, in an interval $M \leq p \leq M + L$ with $M$ and $L$ large such that $L/M$ is a positive
constant smaller than unity.
The simulations strongly suggest  that the structure factor exhibits many well-defined Dirac-delta-function (Bragg-like) peaks along with a small ``diffuse"
contribution; see Fig. \ref{Sofk-numerical}. This means that the primes are characterized by a substantial amount of order
on many length scales, especially relative to the uncorrelated lattice gas (i.e., Poisson distribution of points on the integer lattice) that does not have any such peaks; see Sec. \ref{sec:definitions} for a precise definition.
Motivated by this numerical study,  we  employ analytic number theory to understand rigorously the
nature of the primes as a point process by quantifying the pair correlation function,  structure factor,
local number variance and a certain scalar order metric. While some of the major results were announced
in a letter \cite{To18b}, few mathematical details and derivations were presented there.
Here such details are provided and we also report results that are not contained in Ref. \cite{To18b}.

In Sec. 2, we provide relevant definitions. 
In Sec. 3, we analyze the period-doubling chain,  a simple example of a point process with dense Bragg peaks (Dirac delta function)
that illustrates the phenomenon of \emph{limit-periodicity}, which we will see applies in a modified form to the primes.
In Sec. 4, we {show} that the structure factor defined in (\ref{Skdef1}) is characterized by sharp peaks at certain rational multiples $\pi$, which become progressively denser as $M$ increases, and negligibly small elsewhere. 
The major result is stated in Proposition 1 and a corresponding corollary.
Assuming the Hardy-Littlewood conjecture, this gives what we call an {\it effective limit-periodic} form for the structure factor. 
In Sec.~\ref{sec:hyper}, we show that, in the infinite-size limit, the primes in a dyadic interval form a hyperuniform point process of class II (see Proposition 2). This involves the aforementioned structure factor as well as a cumulative version of it, defined by (\ref{Zofk}), and the number variance
$\sigma^2(R)$ associated with a ``window" of length $2R$. 
In Sec. \ref{sec:tau}, we employ a scalar order metric $\tau$, derived from the structure factor, 
to determine how $\tau$ scales with the system size $L$ (see Proposition 3) and
to identify a transition between large values of $\tau$, when $L$ is comparable to $M$, and small $\tau$ in Gallagher's uncorrelated regime, where $L$ is only logarithmic in $M$. 
In Sec. \ref{sec:class}, we summarize the classification of the primes as a certain limit-periodic, hyperuniform point process.
In Sec. \ref{sec:values}, we describe further numerical investigations into the size of the structure factor $S(k)$.
In Sec. \ref{sec:recon}, we discuss the possibility of reconstructing the primes from the limit-periodic form of the inner sum in Eq.~(\ref{Skdef1}).

\section{Definitions} \label{sec:definitions}

A stochastic point process  in a set $X$ is a collection of points
with configurational positions ${\bf x}_1, {\bf x}_2, {\bf x}_3\ldots$ whose distribution
is described by a probability measure on the set of all possible collections.
Each configuration  in $X$  satisfies two regularity conditions: (i) there are no multiple points
(${\bf x}_i \neq {\bf x}_j$ if $ i\neq j$) and (ii) each bounded subset of $X$ must contain
only a finite number of points. Here we restrict ourselves to one-dimensional
point processes. The set $X$ can be one-dimensional Euclidean space $\mathbb{R}$ (continuum
systems), discrete systems (e.g., the integer  lattice $\mathbb{Z}$), the one-dimensional torus 
$\mathbb{T}$, or discrete systems on  $\mathbb{T}$. The latter two cases constitute {\it periodic}
point processes. A particular configuration (realization) of a point process
in $X$ can formally be characterized by the random variable
\begin{equation}
\n({\bf r})=\sum_{i=1} \delta({\bf r} -{\bf x}_i)
\label{n}
\end{equation}
called the ``local" density at position $\bf r$, where $\delta({\bf r})$
is a $d$-dimensional Dirac delta function. Two particularly important averages
are the one-particle and two-particle correlation functions, $\rho_1({\bf r}_1)$
and $\rho_2({\bf r}_1,{\bf r}_2)$, respectively. When $X$ is $\mathbb{R}^d$ (continuous
systems), they are defined as follows:
\begin{equation}
\rho_1({\bf r}_1)= \langle 
\n({\bf r}_1) \rangle,
\end{equation}
\begin{equation}
\rho_2({\bf r}_1,{\bf r}_2)= \langle \n({\bf r}_1) \n({\bf r}_2)\rangle -\rho_1({\bf r}_1) \delta({\bf r}_1 -{\bf r}_2),
\end{equation}
where the angular brackets denote an average with respect to the probability measure.
The random setting when $X$ is $\mathbb{R}$ is perfectly general and includes 
lattices and periodic point processes as special cases. 
A {\it lattice} in $\mathbb{R}$ is a subgroup
consisting of the integer linear combinations of vectors that constitute a basis for $\mathbb{R}$.
In a lattice in $\mathbb{R}$, the space can be geometrically divided into identical regions called {\it fundamental cells}, each 
of which contains just one point.    In one dimension, there is only one lattice, namely, the integer
lattice $\mathbb{Z}$. The dual of the integer lattice with fundamental-cell spacing $a$ is an integer lattice
with spacing $2\pi/a$, which we denote by $\mathbb{Z}^*$. A one-dimensional periodic point process (crystal) in $\mathbb{R}$ (points in $\mathbb{T}$) is
obtained by placing a fixed configuration of $N$ points (where $N\ge 1$)
within a fundamental cell $F$ of the integer lattice, which
is then periodically replicated.

In the special case of statistically homogeneous point processes in $\mathbb{R}$,
all of the correlation functions are translationally invariant, the first
two of which are then simply given by
\begin{equation}
\rho_1({\bf r}_1)= \rho,
\label{rho1}
\end{equation}
\begin{equation}
\rho_2({\bf r}_1,{\bf r}_2)= \rho^2 g_2({\bf r}_2 -{\bf r}_1).
\label{rho2}
\end{equation}
Here the constant $\rho$ is the number density (number of points per unit volume)
and $g_2({\bf r})$ is the pair correlation function. It is useful to introduce the total correlation function $h({\bf r})$, which
is related to the pair correlation function via
\begin{eqnarray}
h({\bf r})\equiv g_2({\bf r})-1
\label{total}
\end{eqnarray}
and decays to zero for large $|{\bf r}|$ in the absence of long-range order. 
Note that $h({\bf r})=0$ for
all $\bf r$ for the translationally invariant Poisson point process.

The structure factor $S({\bf k})$ is defined as follows: 
\begin{equation}
S({\bf k}) =1+\rho{\tilde h}({\bf k}),
\label{factor}
\end{equation}
where 
\begin{equation}
{\tilde h}({\bf k})= \int_{\mathbb{R}^d} h(\mathbf{r}) \exp\left[-i(\mathbf{k}\cdot  \mathbf{r})\right] d\mathbf{r}
\end{equation}
is the Fourier transform of $h(\bf r)$ so that
\begin{equation}
h({\bf r})= \frac{1}{(2\pi)^d}\int_{\mathbb{R}^d} {\tilde h}(\mathbf{k}) \exp\left[i(\mathbf{k}\cdot  \mathbf{r})\right] d\mathbf{k}.
\label{total-2}
\end{equation}
While $S({\bf k})$ is a nontrivial function for spatially correlated point processes,
it is identically equal to 1 for all $\bf k$
for a translationally invariant Poisson point process. 

In general, the structure factor of a
statistically homogeneous point process can be uniquely be decomposed into three contributions \cite{Ba11}:
\begin{equation}
S({\bf k}) = S({\bf k})_{pp} +  S({\bf k})_{sc} + S({\bf k})_{ac},
\end{equation} 
where $S({\bf k})_{pp}$ is the ``pure point" (Dirac-delta masses) part, 
$S({\bf k})_{sc}$ is the singular-continuous part, and  
$S({\bf k})_{ac}$ is the absolutely-continuous part. In the case of the integer lattice,
$ S({\bf k})$ only consists of the pure-point part. The same is true for a one-dimensional
quasicrystal, such as the Fibonacci chain that is characterized by the golden ratio, except here the Dirac-delta functions are dense \cite{Lev86}.
The Fibonacci chain is a special case of one-dimensional patterns constructed from substitution rules involving
algebraic numbers, and {\it limit-periodic} chains are closely related patterns but are characterized by rational
numbers \cite{Ba11}.  One-dimensional point sets generated from substitution rules involving non-Pisot numbers 
will consist only of singular-continuous contributions \cite{Bom86}.
In the case of a Poisson point process, the only contribution to the structure factor
is the absolutely continuous part. In stark contrast, we will show that the
primes in certain intervals are dominated by a set of dense Bragg peaks.

A hyperuniform statistically homogeneous  point process in $d$-dimensional Euclidean space $\mathbb{R}^d$  
is one in which the structure factor $S({\bf k})$ tends to zero  as the wavenumber $k\equiv |\bf k|$ tends to zero, i.e.,
\begin{equation}
\lim_{|{\bf k}| \rightarrow 0} S({\bf k}) = 0,
\label{hyper}
\end{equation}
implying that single scattering of incident radiation at infinite wavelengths is completely suppressed.
This class of point configurations
includes perfect crystals,  a large class of perfect quasicrystals \cite{Za09,Og17}
and special disordered many-particle systems. Observe that the structure-factor definition (\ref{factor}) 
and the hyperuniformity requirement (\ref{hyper}) dictate that the
volume integral of  $\rho h({\bf r})$ over all space is exactly
equal to $-1$, i.e.,
\begin{equation}
\rho \int_{\mathbb{R}^d} h({\bf r}) d{\bf r}=-1,
\label{sum-1}
\end{equation}
which is a direct-space {\it sum rule} that a hyperuniform point process must obey. 
The hyperuniformity property can be stated in terms of the 
the local number variance $\sigma^2(R)$ associated
within  an interval (window) of length $2R$ for
a one-dimensional homogeneous point process \cite{To03a}:
\begin{eqnarray}
\sigma^2(R)&=& \rho 2R\Big[1+\rho  \int_{\mathbb{R}} h({\bf r}) 
\alpha_2(r;R) d{\bf r}\Big] \nonumber \\
&=&
\frac{\rho R}{\pi} \int_{\mathbb{R}} S({\bf k})
{\tilde \alpha}_2(k;R) d{\bf k},
\label{local}
\end{eqnarray}
where $\alpha(r;R)=1- r/(2R)$ for $r \le 2R$ and zero otherwise, and ${\tilde \alpha}(k;R)= 2 \sin^2(k R)/(kR)$.
A hyperuniform point process is one in which $\sigma^2(R)$
grows more slowly than $R$ in the large-$R$
limit. Three classes of hyperuniformity are to be distinguished: class I, where $\sigma^2(R)$ is bounded; 
class II, where $\sigma^2(R)$ is logarithmic in $R$; and class III, where $\sigma^2(R)$ scales as a power $R^{1-\alpha}$ with $0 < \alpha < 1$ (or $d - \alpha$ for a $d$-dimensional system) \cite{To18a}. After integrating by parts, the second line of (\ref{local}) 
leads to an alternative representation of the number variance \cite{Og17}:
\begin{equation}
\sigma^2(R)=
-\frac{\rho R}{(\pi)} \int_0^\infty  Z(k) 
\frac{\partial {\tilde \alpha}_2(k;R)}{\partial k} dk,
\label{eqn:local-1}
\end{equation}
where
\begin{equation}
Z(K) = 2 \int_0^K  S(k)  dk
\label{Zofk}
\end{equation}
is the {\it integrated} or {\it cumulative} intensity function within a ``sphere" of radius $K$
of the origin in reciprocal space. The quantity  $Z(k)$
has advantages over $S(k)$ in the characterization of quasicrystals
and other point processes with dense Bragg peaks \cite{Og17}.
If $S(k)$ tends to 0 as a power $k^{\alpha}$, then its integral $Z(K)$ will tend to 0 as a power one higher, $Z(K) \sim K^{\alpha+1}$. Any positive power $\alpha > 0$ yields hyperuniformity and distinguishes the system from a random configuration of Poisson points with the same density.

When $X$ is discrete, such as the integer lattice, it sometimes convenient to use the same notation as Eqs.~(\ref{n})-(\ref{rho2})
such that  $\delta({\bf r} -{\bf x}_i)$ is interpreted to be the Kronecker delta
$\delta_{{\bf r},{\bf x}_i}$, which means that $\n({\bf r})$  takes 
either the value 0 or 1, depending on whether the site ${\bf r} \in X$ is unoccupied (empty) or
occupied. In the special case of statistically homogeneous  point processes, while
the definition (\ref{total}) remains the same,
relations (\ref{rho1}) and (\ref{rho2}) are modified  as follows:
\begin{equation}
\rho_1({\bf r}_1)= f,
\end{equation}
\begin{equation}
\rho_2({\bf r}_1,{\bf r}_2)= f^2 g_2({\bf r}_2 -{\bf r}_1),
\end{equation}
where $f$ is the occupation fraction (fraction of occupied sites). Similarly, equation (\ref{factor}) for the structure 
factor becomes in the discrete setting
\begin{equation}
S({\bf k}) =1-f+f{\tilde h}({\bf k}),
\label{factor_d}
\end{equation}
where ${\tilde h}({\bf k})$ is the discrete Fourier transform
\begin{equation}
{\tilde h}({\bf k})= \sum_{\mathbf r\neq \mathbf 0, \mathbf r \in X} h(\mathbf{r}) \exp\left[-i(\mathbf{k}\cdot  \mathbf{r})\right],
\end{equation}
where $h({\bf r})=g_2({\bf r})-1$.
Note that $1-f$ is the structure factor of the {\it uncorrelated lattice gas}, which is a stochastic
point process in $X$ in which the occupation of each site is a constant
probability $f$, independent of any other other site.
While the Fourier-space hyperuniformity condition for discrete $X$ is still given by
relation (\ref{hyper}), the corresponding  direct-space condition
\begin{equation}
 \sum_{\mathbf r\neq \mathbf 0, \mathbf r \in X} h({\bf r})
=\frac{f-1}{f}.
\label{sum-1_dis}
\end{equation}
The relation between the local number variance $\sigma^2(R)$ and the structure factor $S(\mathbf k)$ is unchanged from Eq.~(\ref{local}).
It is noteworthy that the sum rule (\ref{sum-1_dis}) is a condition for hyperuniformity in the
grand-canonical (open-system) ensemble in which the number of particles fluctuates
around some average value; see Refs. \cite{To03a} and \cite{To18a} for details in the 
continuous-space setting. For a system in which the number of particles is fixed, the sum rule still applies
but it is satisfied whether the system is hyperuniform or not.

For  a  single periodic point configuration of $N$ points
within $F$, specified by its local density $\n({\bf r})$ [cf. (\ref{n})], 
it is useful  to introduce
the {\it complex collective density variable} ${\tilde \n}({\bf k })$,
which is simply the Fourier transform of $\n({\bf r})$, i.e.,
\begin{equation}
{\tilde \n}({\bf k }) = \sum_{j=1}^{N} \exp(-i{\bf k \cdot r}_j).
\label{etak}
\end{equation}
This quantity is directly linked to the {\it scattering intensity} ${\cal S}({\bf k})$ defined as
\begin{equation}
{\cal S}({\bf k})= \frac{|{\tilde \n}({\bf k})|^2}{N},
\label{scatter}
\end{equation}
which is a nonnegative real function with inversion-symmetry, i.e., 
\begin{equation}
{\cal S}({\bf k}) = {\cal S}(-{\bf k})
\label{inv}
\end{equation}
that obeys the bounds
\begin{equation}
0 \le {\cal S}({\bf k}) \le N   \qquad ({\bf k} \neq {\bf 0})
\end{equation}
with ${\cal S}({\bf 0})=N $.
For a single periodic configuration with a finite number of
$N$ points within a fundamental cell $F$, the scattering intensity  ${\cal S}({\bf k})$
is identical to the structure factor $S({\bf k})$ [cf. (\ref{factor})], except the latter excludes $\bf k=0$ (forward
scattering). In general, whether they remain equal in the infinite-system limit depends on the ergodicity of the process, but this issue does not affect our analysis of the primes, and so we will simply take Eq.(\ref{Skdef1}) to be the definition of the structure factor.
Importantly, the definition of hyperuniformity excludes the forward scattering contribution, which is implicit in (\ref{hyper}).


A useful scalar positive order metric that is capable of capturing the 
degree of translational order across length scales is the $\tau$ order metric \cite{To15}. For  a statistically
homogeneous point process in $\mathbb{R}^d$ at number density $\rho$, it is defined by
\begin{eqnarray}
\tau &\equiv &\frac{1}{D^d} \int_{\mathbb{R}^d} [g_2({\bf r})-1]^2 d\mathbf r \\
&=& \frac{1}{(2\pi)^d D^d}\int_{\mathbb{R}^d} [S({\bf k})-1]^2 d\mathbf k,
\label{tau_c}
\end{eqnarray}
where $D$ is some characteristic length scale.
A convenient choice is $D=\rho^{-1/d}$. For a Poisson point process in $\mathbb{R}^d$, $\tau=0$ 
because $g_2({\bf r})-1$ is zero  for all $\bf r$. Thus, a deviation of $\tau$ from
zero measures translational order with respect to the fully uncorrelated case.
For example, for the Riemann zeta zeros, $\tau=2/3$,  assuming Montgomery's pair correlation conjecture or, equivalently,
the corresponding structure factor (\ref{spec}), which reflects the disordered hyperuniformity of the point process.
For any periodic point process in which there are a finite number of points
within the fundamental cell $F$, $\tau$ is unbounded because the integrals
are carried out over all space. For this reason, one can employ a modified
version of $\tau$ by carrying out the integral in direct space
or reciprocal space over appropriate subsets
of $\mathbb{R}^d$, in which case the equality in Eq. (\ref{tau_c})
no longer applies. 


The discrete-setting counterpart of the order metric $\tau$ defined in (\ref{tau_c}) in which $X$ 
is a subset of $\mathbb{Z}$ on the torus $\mathbb{T}$ in which the fundamental cell has length $L$ is given by 
\begin{eqnarray}
\tau &\equiv& \sum_{j=1}^{Ns-1} f^2 [g_2(2j)-1]^2\\
&=& \frac{1}{N_s} \sum_{j=1}^{Ns-1}\left(S\left(\frac{j\pi}{N_s}\right)-(1-f)\right)^2,
\label{tau-discrete}
\end{eqnarray}
where $N_s$ is the number of lattice sites within the fundamental cell and $N$ is the number of occupied sites.
Strictly speaking, the quantities $g_2(r)$  and $S(k)$ are ensemble averages.
In the case of an uncorrelated lattice gas, $S=1-f$  in the infinite-system-size limit
so that $\tau=0$.
The corresponding expression in the discrete setting for a single configuration
in any space dimension was presented and applied in Ref. \cite{Di18}. Here $g_2(r)$  and $S(k)$
should be interpreted to come from a single configuration. Analysis of the primes in some fixed interval requires
the use of this single-configuration variant of $\tau$, which we will employ. Note that $\tau$
for a single configuration of an uncorrelated lattice gas is given by $(1-f)^2$ (not zero) in the infinite-system-size limit.


\section{An Illustrative Example: The Period-Doubling Chain}
\label{doubling}

The spatial distribution of the primes shares some features with \emph{limit-periodic} point sets, so we discuss a model example in detail. Consider two types of intervals (``tiles" or ``letters"): $a$ and $b$. 
The period-doubling chain is defined by the following iterative substitution rule: $a \rightarrow ab$ and $b \rightarrow aa$ \cite{Ba11}.
In the infinite-size limit, this constitutes a point process on $\mathbb{Z}$ in which a subset of sites are occupied by $a$'s
and the remaining sites are occupied by $b$'s.
The locations of the $b$'s are given by a superposition of arithmetic progressions $2+4j$, $8+16j$, $32+64j$, with a factor of 4 from one to the next. Thus the infinite-size limit is a union of periodic systems, which is
termed limit-periodic. The limiting densities of $a$ and $b$ sites are 2/3 and 1/3, respectively. The structure factor associated with the $a$'s is given by
\begin{equation}
S(k) = \frac{4\pi}{3}\left(\sum_{m=1}^{\infty} \delta(k-2\pi m) + \sum_{n=1}^{\infty} \sum_{m=1}^{\infty} 2^{-2n} \delta\left(k-\frac{(2m-1)\pi}{2^{n-1}}\right) \right)\,
\label{S-doubling}
\end{equation}
assuming unit lattice spacing. This is obtained by squaring Eq.~(12) in Ref. \cite{Ba11}, multiplying it by $2\pi/f$, and then rescaling the function by $2\pi$.  The factor $2\pi$ accounts for differences in the definition of the Fourier transform. Thus, we have a dense set of Dirac-delta-function peaks, one for each dyadic rational $(2m-1)/2^{n-1}$; see Fig. \ref{S-pd}. These peaks at certain rational numbers arbitrarily close to 0 are a feature shared by this example and the prime numbers. Figure \ref{S-pd} depicts the structure factor of the period-doubling chain.

\begin{figure}[bthp]
\includegraphics[  width=3in,keepaspectratio,clip=]{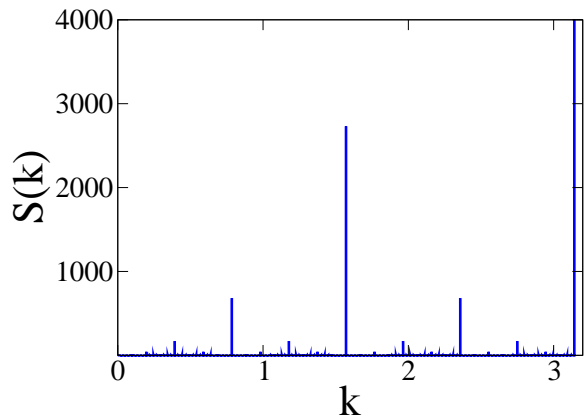}
\caption{Structure factor of the period-doubling chain as obtained from formula (\ref{S-pd}) with $n=20$. Note the self-similarity in the intensities
and locations of the peaks. The height of the peak at $k=\pi$ is larger than shown on the scale of this figure. }
\label{S-pd}
\end{figure}

Substitution of this expression ({\ref{S-doubling}}) for the structure factor into relation ({\ref{local})
yields the local number variance for the period-doubling chain:
\begin{equation}
\sigma^2(R) = \frac{8}{9\pi^2}\left(\sum_{m=1}^{\infty} \frac{\sin^2(2m\pi R)}{m^2} + \sum_{n=1}^{\infty}\sum_{m=1}^{\infty} \frac{\sin^2((2 m-1)\pi R/2^{n-1})}{(2m-1)^2} \right).
\label{pd}
\end{equation}
The first term in (\ref{pd}) is a periodic function $R(1-2R)/9$ with period $1/2$ 
and the second term  
is a superposition of periodized ``triangle" functions with heights $1/9$ and bases $1$, $2$, $4$, $\cdots$..
Together this results in a number variance that grows logarithmically in $R$; see Fig. \ref{sigma-pd}.
Therefore, the period-doubling chain is hyperuniform of class II.

\begin{figure}[bthp]
\includegraphics[  width=3in,keepaspectratio,clip=]{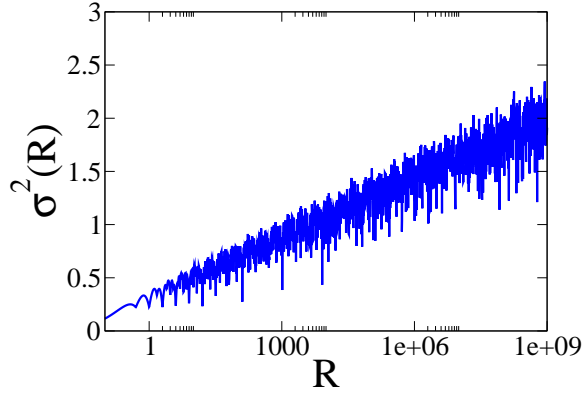}
\caption{Local number variance of the period-doubling chain as computed from the explicit formula (\ref{pd}) with $n=20$.}
\label{sigma-pd}
\end{figure}

Having established that the period-doubling chain is hyperuniform, we now want to determine how the structure factor $S(k)$ behaves in the vicinity of the origin.
The structure factor $S(k)$ is not a continuous function because there are dense Bragg peaks arbitrarily close to 0, so we do not have $S(k) \rightarrow 0$ as $k \rightarrow 0$ in the usual sense. We follow the practice of Ref. \cite{Og17} in such instances and pass to a cumulative version of the structure factor,
$Z(K)$, defined by (\ref{Zofk}).
This integral simply adds the weights of the $\delta$ peaks up to position $K$. In order to have a peak $(2m-1)\pi/2^{n-1}$ within the range of integration, $n$ must be relatively large:
\begin{equation}
\exists m \ \frac{2m-1}{2^{n-1}} \pi < K \Longleftrightarrow 2^{n-1} > \pi/K \Longleftrightarrow n > \log_{2}(\pi/K) + 1 = \log_{2}(2\pi/K),
\end{equation}
or else there are no integers $m$ in the necessary interval. We have an explicit formula for the tail of a geometric series:
\begin{equation}
\sum_{n > C} b^n = b^{\lceil C \rceil} \frac{1}{1-b}
\end{equation}
where $\lceil C \rceil$ denotes the least integer greater than $C$ (and, in particular, $C+1$ in case $C$ is already an integer).
Using this to sum the series in $Z(K)$ leads to an explicit formula
\begin{align}
\label{Z-pd}
Z(K) &= 2\int_0^K S(k)dk = \frac{8\pi}{3} \sum_{n > \log_2 (2\pi/K)} 2^{-2n} \left\lfloor \frac{1}{2} + 2^{n-2}\frac{K}{\pi} \right\rfloor \\
&= \frac{8\pi}{3} \sum_{n > \log_2 (2\pi/K)} 2^{-2n} \left( \frac{1}{2} + 2^{n-2} \frac{K}{\pi} + \left\{\frac{1}{2} + 2^{n-2}\frac{K}{\pi} \right\}\right) \nonumber \\
&= \frac{8\pi}{3} \left( \frac{4}{3} 2^{-2 \lceil \log_2(2\pi/K) \rceil} + \frac{1}{2}\frac{K}{\pi} 2^{-\lceil \log_2(2\pi/K) \rceil} + \sum_{n > \log_2(2\pi/K)} 2^{-2n} \left\{ \frac{1}{2} + 2^{n-2}\frac{K}{\pi} \right\} \right) \nonumber
\end{align}
where the braces $\{ \cdot \}$ denote fractional part. Taking into account the jump discontinuities 
when $K/\pi$ crosses a dyadic rational shows that
\begin{equation}
\frac{1}{6\pi}K^2 \leq Z(K) \leq \frac{1}{2\pi} K^2.
\label{Z-UL}
\end{equation}
Figure \ref{Z-pd-plot} shows the function $Z(k)$ and the aforementioned upper and lower bounds.


Thus $Z(K)$ is bounded between two multiples of $K^2$, with a self-similar staircase-like behaviour in between. 
Using these bounds and relation (\ref{eqn:local-1}), we get the following corresponding asymptotic bounds on the number variance $\sigma^2(R)$ in the limit $R \rightarrow \infty$:
\begin{equation}
\frac{4}{9\pi^2} \ln(R) \le \sigma^2(R) \le \frac{4}{3\pi^2} \ln(R).
\end{equation}
This implies that the period-doubling chain falls within class II of
hyperuniform systems with a structure factor that effectively behaves as $S(k) \sim k$ as $k \to 0$ \cite{To18a}.}
These asymptotic bounds closely match the upper and lower envelopes of the fluctuating number variance function plotted in Fig. \ref{sigma-pd}.

\begin{figure}[bthp]
\includegraphics[  width=3in,keepaspectratio,clip=]{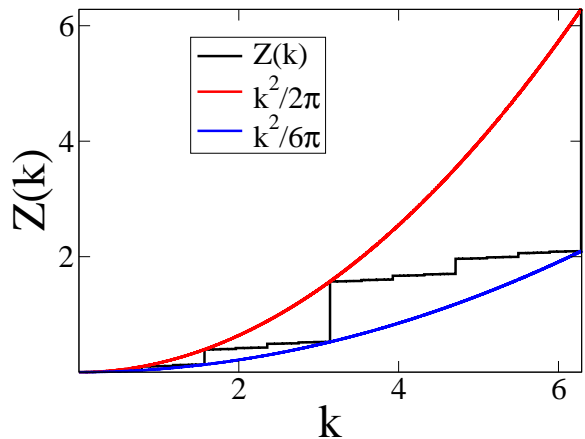}
\caption{The cumulative intensity function $Z(k)$ of the period-doubling chain as obtained from formula (\ref{Z-pd}) with $n=20$. The upper and lower bounds given in (\ref{Z-UL}) are also indicated in the figure.}
\label{Z-pd-plot}
\end{figure}

\section{Primes in Progressions and Peaks in the Structure Factor} \label{sec:circ}

Here we provide the theoretical basis for the numerical results reported in Ref. \cite{Zh18}
that shows that the dominant contribution to the structure factor $S(k)$ 
in an interval $[M,M + L]$ (with $M$ and $L$ large, and $L/M$ smaller than unity) consists of many 
well-defined Dirac-delta-function peaks. We show that in the infinite-system-size limit 
($M \to \infty$, $L/M\to \beta >0$) and after scaling by the density $\rho$, 
the peaks in the structure factor of the primes will become Dirac delta functions
at rational numbers with odd, square-free denominators, and hence the small diffuse part
observed numerically in  Ref. \cite{Zh18}, vanishes in this limit. This major result is summarized
in the following proposition:
\bigskip
 
\noindent
{{\sl Proposition 1:} The structure factor of the prime numbers, scaled by the 
the density $\rho$, in an interval $M \leq p \leq M + L$ 
in  the limit such that $M \rightarrow \infty$, $L/M \rightarrow \beta > 0$ is given by
\begin{equation}
\lim_{M \rightarrow \infty} \frac{S(k)}{2\pi \rho} = {\sum_{n}}^{\flat} {\sum_{m}}^{\times} \frac{1}{\phi(n)^2} \delta\left(k - \frac{m\pi}{n}\right).
\label{S-p}
\end{equation}
Here, $\phi(n)$ is Euler's totient, which counts the numbers up to $n$ with no factor in common with $n$ \cite{Te95},
the symbol $\flat$ indicates that the sum over $n$ only involves odd, square-free values of $n$ 
and the symbol $\times$ indicates that $m$ and $n$ have no common factor.}
\bigskip

\noindent
{{\it Remark 1:} 
Since $L/M=\beta$ is fixed, and in particular $M+L$ is within a constant multiple of $M$, the prime number theorem implies that 
the density $\rho$ is effectively $1/\ln(M)$ throughout the interval.
We see  that this normalized structure factor of the primes (\ref{S-p}) exhibits dense Bragg peaks 
at certain rational values of $k/\pi$, and hence is similar to the structure factor (\ref{S-doubling}) of the limit-periodic period-doubling chain 
discussed in Sec. \ref{doubling}. However, there is a fundamental difference between these two systems, to be elaborated on in Sec. \ref{class}.}

\bigskip

\noindent
{{\it Remark 2:} The interpretation of (\ref{S-p}) is that $S(k)/(2\pi \rho)$ converges to the sum of peaks $\sum_n \sum_m \phi(n)^{-2} \delta(k - \pi m/n)$ in the sense that their integrals against test functions $f(k)$ are close:}
\begin{equation} \label{S-p-test}
\frac{1}{2\pi \rho } \int_0^{\pi} f(k)  S(k) dk \approx {\sum_{n}}^{\flat} {\sum_{m}}^{\times} \frac{1}{\phi(n)^2} f(\pi m/n)
\end{equation}
In particular, consider $f(k) = e^{-irk}$. For $r \neq 0$, this gives a count of how often $p$ and $p+r$ are both prime. 
\bigskip

{In what follows, we give a derivation of (\ref{S-p}). Our approach
is somewhat different than the original approach of Hardy and Littlewood,
which is outlined in Appendix A.} 


\subsection{{Derivation of Proposition 1}}
Our first step is to replace $S(k)$ by another sum involving the more convenient weights $\Lambda(n)$ given by $\ln(p)$ when $n$ is a power of the prime $p$, and 0 otherwise. One uses summation by parts to convert the weights, and then includes higher powers at essentially no cost.
Indeed, if $p^t < Y$ with $t\geq 2$, then $p < Y^{1/2}$. A sum over so few terms introduces an error of order no worse than $L^{1/2}$. 
{The result
for the complex density variable ${\tilde \n}(k)$ (cf. \ref{etak}) is}
\begin{equation}
{\tilde \n}(k)=\sum_{M<p\leq M+L} e^{ik p} = \frac{1}{\ln(M)} \sum_{n=M+1}^{M+L} \Lambda(n) e^{i\alpha n} + O\left(\frac{L}{\ln(M)^2} \right).
\end{equation}
Let us take the absolute square and, to bound the error, resort to the trivial bound $\sum \Lambda(n) e^{i\alpha n} \lesssim L$. We divide by $\rho N$ and note that $\rho \sim 1/\ln(M)$ and $N \sim \rho L$, to arrive at
\[
\frac{1}{2\pi \rho} S(k) = \frac{1}{2\pi L} \left| \sum_{n=M+1}^{M+L} \Lambda(n) e^{ikn} \right|^2 + O\big(L/\ln(M)^2 \big).
\]


We split the integral into arcs $\mathfrak{M}(q,a) = \{ |\alpha - a/q | < \varepsilon \}$ near fractions $a/q$ with small denominator $q$ (major arcs), and write $\mathfrak{m}$ for the rest of the interval (minor arcs). The denominator is restricted to $q \leq q_{\max}$ and we choose $q_{\max} = \ln(L)^B$, with a constant $B$ as large as one pleases. The length $\varepsilon$ of each major arc is chosen to be $\varepsilon = \ln(L)^B/L$, and in principle different lengths could be adapted to the test function $f$. Let us first calculate the contribution from the arcs $\mathfrak{M}(q,a)$, returning later to the discussion of the minor arcs. 

On the major arc $\mathfrak{M}(q,a)$, even without the Riemann hypothesis, we still have the approximation
\begin{equation} \label{eq:trig}
\sum_{n=M+1}^{M+L} \Lambda(n) e^{ikn} = \frac{\mu(q)}{\phi(q)} \sum_{n=M+1}^{M+L} e^{i(k - 2\pi a/q)n} + O\big(Le^{-c \sqrt{L}} \big)
\end{equation}
This plays the role of Hardy-Littlewood's approximation to $f(x)$ in Lemma 9 of \cite{Ha23}.
It follows from p. 147 in \cite{Da80} or, what is the same, Lemma 3.1 in \cite{vaughan97} (p.30). They take $M=1$ but one can of course subtract. We summarize the argument from \cite{Da80} for the reader's benefit in order to emphasize the role played by primes in arithmetic progression. 

The basic idea is to decompose with respect to Dirichlet characters modulo $q$. These are functions $\chi(n)$ defined for integer values of $n$ and characterized by the properties of periodicity and multiplicativity, namely $\chi(n + q) = \chi(n)$ and $\chi(ab)=\chi(a)\chi(b)$. They are the natural harmonics to use in a situation with period $q$ that preserves multiplicative structure, such as remainders after division by $q$. One character $\chi_0$ is distinguished, given by $\chi_0(n) = 1$ for $\gcd(n,q)=1$ and $\chi_0(n)=0$ in case $\gcd(n,q)>1$. 
As a rule of thumb, $\chi_0$ provides the main term and we must endeavor to show that the contributions from other characters $\chi$ are negligible. 
The phase $e(an/q)$ can be expanded as a sum over Dirichlet characters $\chi$ modulo $q$. Assuming $\gcd(an,q)=1$, we have
 \begin{equation} \label{eq:characters}
 e(an/q) = \frac{1}{\phi(q)} \sum_{\chi} G(\overline{\chi}) \chi(ap).
 \end{equation}
The sum is over all Dirichlet characters modulo $q$. If $\gcd(an,q) > 1$, then the sum is 0.
We recommend \cite{Da80} for the theory of Dirichlet characters as well as the Gauss sum $G(\overline{\chi})$. By definition, 
 \begin{equation}
 G(\chi) = \sum_{m=1}^q \chi(m) e^{2\pi i m/q},
 \end{equation}
and it is important to note that $|G(\chi)| \leq q^{1/2}$ while $G(\chi_0) = \mu(q)$ is the M\"{o}bius function.
Since $n$ is a prime power and $a$ has no factor in common with $q$, $\gcd(an,q)=1$ does hold unless $n$ is a power of a prime dividing $q$. Factoring $q$ as $q = p_1 p_2 \cdots p_w \geq 2^w$ shows that there are at most $\ln(q)$ primes dividing $q$. For each such $p$, a prime power $n=p^l$ will be less than $M+L$ for $l \lesssim \ln(M)$. In our range, with a value of $q_{\max}$ much less than $M$, we thus have $\gcd(an,q)=1$ except for at most $\ln(M)^2$ terms $n$. With an error no worse than a power of $\ln(M)$, we may thus ignore these exceptions, proceeding as if (\ref{eq:characters}) held for all $n$.
Write $\alpha = k/\pi = a/q + t$ and $e(z) = e^{2\pi i z}$ for convenience. Summing over $n$ gives
 \begin{equation}
 \sum_{n=M+1}^{M+L} e(\alpha n) \Lambda(n) = \frac{1}{\phi(q)} \sum_{\chi} \tau(\overline{\chi}) \chi(a) \sum_{n=M+1}^{M+L} e(nt)\chi(n)\Lambda(n).
 \end{equation}
 Summing by parts, we have
 \begin{equation}
 \sum_{n \leq Y} e(nt)\chi(n)\Lambda(n) = e(Yt) \psi(Y,\chi) - 2\pi i t \int_1^Y e(ut)\psi(u,\chi) du,
 \end{equation}
 where, for an upper limit $u$ and a character $\chi$ to modulus $q$,
 \begin{equation}
 \psi(u,\chi) = \sum_{n < u} \Lambda(n)\chi(n). 
 \end{equation}
By the prime number theorem in progressions (p. 125, 132 of \cite{Da80}), there is a positive $c>0$ such that for all characters to moduli $q \leq \ln(L)^B$, $\psi(u,\chi_0) = u + O\big(\exp(-c \sqrt{\ln(u)}\big)$ while for nontrivial $\chi$, $\psi(u,\chi) = O\big(\exp(-c\sqrt{\ln(u)})\big)$. This leads to the error term stated in (\ref{eq:trig}) (and the Riemann hypothesis would imply an estimate for an even larger range of $q$).

Summing the geometric series, we have
\begin{equation}
\left|\sum_{n=M+1}^{M+L} e(nt) \right|^2 = \frac{1-\cos(2\pi tL)}{1 -\cos(2\pi t)} = \left(\frac{\sin(\pi Lt)}{\sin(\pi t)}\right)^2.
\end{equation}
Note that $e^{c \sqrt{\ln(L)} }$ is much larger than $\ln(L)^b = e^{b \ln\ln{L} }$ for any positive values of $a$ and $b$. Thus the errors from using the prime number theorem are smaller than the error $L/\ln(M)^2$ that we have already exposed ourselves 
to through summation by parts. Therefore, the structure factor in the vicinity of a particular
peak location $k$ for some $a$ and $q$ and sufficiently small $t$ is given by
\begin{equation}
\frac{1}{2\pi \rho} S(k) = \frac{\mu(q)^2}{\phi(q)^2} \frac{1}{2\pi L} \left(\frac{\sin(\pi Lt)}{\sin(\pi t)}\right)^2 + O\big( L /\ln(M)^2 \big).
\label{S-sinc}
\end{equation}
In the limit that $t$ tends to zero faster than $L$ goes to infinity, and 
$M \to \infty$ such that $L/M$ remains finite,
this formula with $q=2n$ and $a=m$ yields the limiting form (\ref{S-p}) of {\sl Proposition 1}
with dense Bragg peaks.
\bigskip

\noindent
{\sl Corollary:} For finite but large $N$, 
the structure factor at a rational value of $k/\pi$ is given by
\begin{equation}
S(\pi m/n) \sim \frac{N}{\phi(2n)^2} \mu(2n)^2.
\label{SkPeak}
\end{equation}
This formula follows immediately from (\ref{S-sinc}). The M\"{o}bius function $\mu(n)$ \cite{Te95} appears here
because $\mu^2(2n)$
is one whenever $2n$ is square-free and zero otherwise. Also, $\phi(n) = \phi(2n)$ if $n$ is odd.
We will show below that this formula is in excellent agreement with
the numerical results reported in Ref. \cite{Zh18}.

Now we provide strong arguments to demonstrate that the structure factor is negligibly small
at the irrationals in the infinite-system-size limit. We recall Fej\'{e}r's theorem that for a continous, periodic function $f$, we have uniform convergence
\begin{equation}
\frac{1}{2\pi} \int_{-\pi}^{\pi} f(x-t) \frac{1}{L} \left( \frac{\sin(Lt/2)}{\sin(t/2)}\right)^2 dt \rightarrow f(x)
\end{equation}
as $L \rightarrow \infty$. \\
Assuming $\varepsilon$ is not too small, the bulk of the integral comes from $|t| < \varepsilon$. Indeed,
\[
\int_{\varepsilon}^{1-\varepsilon} f(a/q+t) \frac{1}{L} \left( \frac{\sin(\pi Lt)}{\sin(\pi t)} \right)^2 dt \leq \frac{2 \| f \|_{\infty} }{L} \int_{\varepsilon}^{1/2} \frac{dt}{\sin(\pi t)^2} = \frac{2 \|f\|_{\infty}}{L\pi}\cot(\pi \varepsilon) \lesssim \frac{\| f \|_{\infty} }{L\varepsilon}
\]
This implies that the integral over $\mathfrak{M}(q,a)$ converges to $f(a/q)$ as $L\rightarrow \infty$, provided $L\varepsilon \rightarrow \infty$. This is the case, for example, with our choice $\varepsilon = \log(L)^B/L$ for any power $B$. Note that
\[
\frac{1}{L} \int_0^1 \left(\frac{\sin(\pi L t)}{\sin(\pi t)} \right)^2 dt = 1.
\]
This allows us to write
\begin{align*}
& \int_0^1 f(a/q + t) \frac{1}{L} \left(\frac{\sin(\pi L t)}{\sin(\pi t)} \right)^2 dt = f(a/q) + \int_0^1 \big(f(a/q+t)-f(a/q) \big) \left(\frac{\sin(\pi L t)}{\sin(\pi t)} \right)^2 \frac{1}{L} dt \\
&= f(a/q) + \int_{-\varepsilon}^{\varepsilon} \big(f(a/q+t)-f(a/q) \big) \left(\frac{\sin(\pi L t)}{\sin(\pi t)} \right)^2 \frac{1}{L} dt + O\left( \frac{\| f \|_{\infty} }{L\varepsilon} \right)
\end{align*}
Suppose that $f$ has modulus of continuity $\omega$, meaning that $| f(x) - f(y) | \leq \omega(\varepsilon)$ whenever $|x-y| \leq \varepsilon$. For example, if $f$ has continuous derivatives of order up to $m$, then $\omega(\varepsilon) = \| f \|_{C^m} \varepsilon^m$ is a valid modulus of continuity for $f$. Using once more the fact that the integral of the Fej\'{e}r kernel is 1, we have
\begin{align*}
\int_0^1 f(a/q + t) \frac{1}{L} \left(\frac{\sin(\pi L t)}{\sin(\pi t)} \right)^2 dt &= f(a/q) + O\left(\omega(\varepsilon) + \frac{\| f \|_{\infty}}{L\varepsilon} \right) 
\end{align*}
Dividing the interval $[0,1]$ into the major arcs $\mathfrak{M}(a,q)$ for $q \leq q_{\max}$ together with the remaining minor arcs $\mathfrak{m}$, and summing the error made in Fej\'{e}r's theorem over $q \leq q_{\max}$, we get
\begin{align*}
 &\int_0^1 f(\alpha) \frac{1}{L} \left| \sum_{X < n < Y} e(n\alpha)\Lambda(n) \right|^2 d\alpha = \\
 &\sum_{q} {\sum_a}^{\times} \frac{\mu(q)^2}{\phi(q)^2} f\big(\frac{a}{q}\big) + O\left( \big( \omega_f(\varepsilon) + \frac{\| f \|_{\infty} }{L\varepsilon} \big) \ln(q_{\max}) \right) + \frac{1}{L} \int_{\mathfrak{m}} f(\alpha) \left|\sum_{X <n <Y} e(n\alpha)\Lambda(n) \right|^2 d\alpha
\end{align*}
Note  that $\mu(q)^2$} is simply 1 or 0 according to whether $q$ is squarefree, so the sum is restricted to squarefree $q$.
This agrees with the sum of $\delta$ functions in our limiting form (\ref{S-p}). 
We assume that $f$ is smooth enough to have $\omega_f(\varepsilon) \ln(q_{\max}) \rightarrow 0$, where $\epsilon = \ln(L)^B/L = q_{\max}/L$. For example, this is the case if $f$ is $C^1$ or even just H\"{o}lder continuous with any exponent $\alpha$, since $\varepsilon^{\alpha} \ln(q_{\max}) \rightarrow 0$. In particular, one may take $f(k) = e^{ikr}$ for the application to prime pairs. Likewise, $\ln(q_{\max})/L\varepsilon \rightarrow 0$.
Therefore, assuming the integral over $\mathfrak{m}$ is negligible, we expect that
\begin{equation}
 \int_0^1 f(\alpha) \frac{1}{L} \left| \sum_{n=M+1}^{M+L} e(n\alpha)\Lambda(n) \right|^2 d\alpha \approx \sum_q \sum_a \frac{\mu(q)^2}{\phi(q)^2} f(a/q). 
\end{equation}
Let us reiterate that proving that the minor arcs contribute less than the major arcs is a significant challenge, which we have not solved. Nevertheless, the analysis is very suggestive, and the limiting form (\ref{S-p}) is also supported by numerical results, as we now illustrate.

\begin{figure}[bthp]
\includegraphics[  width=4in,keepaspectratio,clip=]{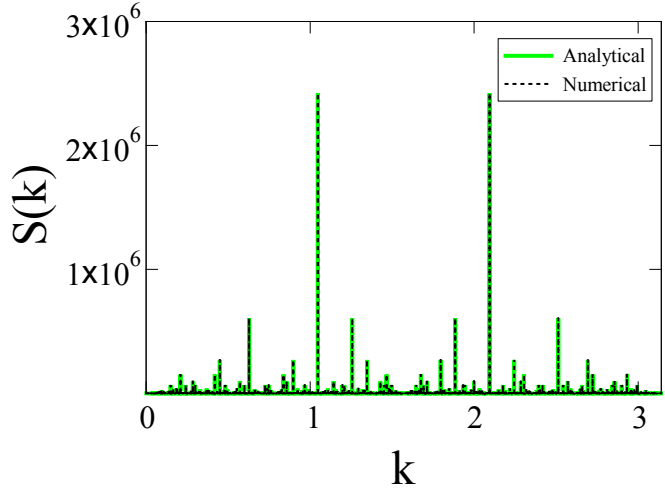}
\caption{Structure factor for the primes  as a function of $k$ (in units of the integer lattice
spacing) in the interval $[M, M+L)$, as predicted from formula (\ref{SkPeak}) for $M=10^{10}+1$, $L=2.23\times10^{8}$ and $n_{max} =100 \ln(M)$.
This shows many peaks with a type of self-similarity described in the main body of the text and is seen to be  in excellent agreement with the corresponding numerically computed structure factor obtained in Ref. \cite{Zh18}.}
\label{S-primes}
\end{figure}

Figure \ref{S-primes} compares the prediction of formula (\ref{SkPeak}) for the structure factor
of the primes for $M=10^{10}+1$, $L=2.23\times10^{8}$ and $n_{max} =100 \ln(M)$
to the corresponding numerical  results reported in Ref. \cite{Zh18}. Agreement between the
analytical and numerical results is excellent. 
The structure factor contains many well-defined Bragg-like peaks of various intensities
characterized by a type of self-similarity.
This is explained by the fact that $\phi(n_1 n_2) = \phi(n_1)\phi(n_2)$ for relatively prime $n_1$ and $n_2$, so that rescaling preserves the relative heights of the peaks given by Eq. (\ref{S-p}). 

\subsection{Pair Correlation Function}
\smallskip

We now obtain the pair correlation function $g_2(r)$ of the primes  via the limiting form
of the structure factor (\ref{S-p})
by performing the inverse Fourier transform of $S(k)-1 \equiv \rho {\tilde h}(k)$ using (\ref{total-2}),
where ${\tilde h}(k)$ is the Fourier transform of $h(r)\equiv g_2(r)-1$ for $r\neq 0$:
\begin{equation}
g_2(r)=1+{\sum_{n}}^{\flat} \frac{1}{\phi^2(n)}{\sum_{m}}^{\times} \exp(rm\pi i/n).
\label{G2}
\end{equation}
Recall that the pair correlation function is defined such that $fg_2(r)$ gives the conditional probability that, assuming $p$ is prime, so is $p+r$.
Because the density of the primes is $1/\ln(M)$, the occupation fraction is $f=2/\ln(M)$, and the density of prime pairs with separation $r$ [cf. (\ref{rho2})] is
\begin{equation}
\begin{split}
\rho(r)&=\frac{\# \{ p \in [M,M+L] \ ; \ p, \ p+r \ \text{prime} \}}{L}\\ 
&=\frac{f^2 g_2(r)}{2} \\
&= \frac{2}{(\ln M)^2}\left[   1+{\sum_{n}}^{\flat} \frac{1}{\phi^2(n)}{\sum_{m}}^{\times} \exp(rm\pi i/n)\right].
\end{split}
\label{pair}
\end{equation}

From (\ref{eq:tuples}) and (\ref{pair}), we see that expression (\ref{G2}) for $g_2(r)$  for distinct values of
$r = 2, 4, 6,\ldots$  is simply a different representation of the 
Hardy-Littlewood constant $\mathfrak{S}(\mathcal{H})$, given by (\ref{HLConstant}),  for the case $m=2$ so that $\mathcal{H} = \{ 0, r \}$.
Thus, (\ref{S-p}) implies the Hardy-Littlewood conjecture on prime pairs by taking the test function $f(k)$
of (\ref{S-p-test}) to be $e^{ikr}$. Conversely, if the Hardy-Littlewood conjecture holds for every $r$, then one knows (\ref{S-p}) for exponential functions $e^{irk}$ and one can deduce it for other functions $f(k)$ by Fourier expansion. The limiting form for $S(k)$ is thus an equivalent formulation of the Hardy-Littlewood conjecture. As such, we certainly do not have a proof of it to offer here.



\section{Hyperuniformity} \label{sec:hyper}
 

\noindent{\sl Proposition 2:} The primes  in the infinite-system-size limit 
such that $M \rightarrow \infty$, $L/M \rightarrow \beta > 0$ form a hyperuniform point process of class II.
\bigskip

\noindent{\it Proof:} As in the period-doubling chain, to determine whether the primes are hyperuniform, we 
first determine the  cumulative intensity $Z(K)$ defined by (\ref{Zofk}). {Using this relation
and (\ref{S-p}), we  find that this quantity
\begin{equation}
\lim_{M \rightarrow \infty} \frac{Z(K)}{2\pi \rho} = 2\,{\sum_{n}}^{\flat} {\sum_{\frac{\pi m}{n} < K}}^{\times} \frac{1}{\phi(n)^2}.
\label{Z-K}
\end{equation}}
Given a cutoff so that only values $cN$ count as peaks in the finite-system structure factor, the denominator $n$ does not exceed $\sqrt{1/c}$. On the other hand, $n$ must be large enough for there to be a peak $m\pi/n$ less than $K$:
\begin{equation}
\frac{\pi}{K} \leq n \leq \sqrt{\frac{1}{c}} = n_{\text{max}}.
\end{equation}
In particular, the lowest allowable $K$ is $\pi/n_{\max}$.

For the sum over $m$, note that
\[
{\sum_{m < cn} }^{\times} 1 = c\phi(n) + O_{\epsilon}(n^{\epsilon}) ,
\]
Thus the number of $m$ coprime to $n$ is growing regularly: Summing only up to $cn$ yields a fraction $c$ of the total $\phi(n)$. To see this, we use the M\"{o}bius inversion property, namely
\begin{equation}
\sum_{d | t} \mu(d) = \begin{cases}
                       1 \ \text{if} \ t = 1 \\
0 \ \text{if} \ t > 1
                      \end{cases}
\end{equation}
to detect the condition $\gcd(m,n)=1$.
This implies that for any $c$ (it will be $K/\pi$ for us, no longer the same as the peak cutoff $c$ above)
\begin{equation*}
{\sum_{m < cn }}^{\times} 1 = \sum_{m < cn} \sum_{d | \gcd(m,n)} \mu(d) 
= \sum_{d | n} \sum_{b < cn/d} \mu(d) 
= \sum_{d | n} \mu(d) \lfloor cn/d \rfloor.
\end{equation*}
Now we write $\lfloor cn/d \rfloor = cn/d - \{ cn/d \}$ in terms of integer part and fractional part:
\[
\sum_{d | n} \mu(d) \lfloor cn/d \rfloor 
= c \sum_{d | n} \mu(d) \frac{n}{d} - \sum_{d | n} \mu(d) \{cn/d\}.
\]
The main term $c\phi(n)$ thus comes from the identity
\begin{equation*}
\phi(n) = \sum_{d | n} \mu(d) \frac{n}{d}.
\end{equation*}
The rest is negligible because $| \mu(d) \{ cn/d \}| \leq 1$ and the number of divisors of $n$ is $O_{\epsilon}(n^{\epsilon})$.

The consequence is that the cumulative structure factor (\ref{Z-K}) becomes 
\begin{equation}
\label{ZkPeak}
\lim_{M \rightarrow \infty} \frac{Z(K)}{2\pi \rho}  = 2\,{\sum_{n}}^{\flat} \frac{1}{\phi(n)^2} \left( \frac{K}{\pi} \phi(n) + O_{\epsilon}(n^{\epsilon})\right)
\end{equation}
To a good approximation for sufficiently large $M$,  this yields 
\begin{align}
Z(K) &\approx \frac{4\,K}{\ln(M)} {\sum_{\frac{\pi}{K} < n < n_{\text{max}}}}^{\hspace{-0.18in}\flat} \frac{1}{\phi(n)} + O_{\epsilon}\left(\sum_{n} \frac{n^{\epsilon}}{\phi(n)^2}\right)\\
&\approx \frac{K}{\ln(M)}\left( \ln{n_{\text{max}}} - \ln{\frac{\pi}{K}}\right) 
\label{Z-primes}
\end{align}
Recall that there is a lower bound on $K$: $K \geq \pi/n_{\max}$. Expanding 
the function $Z(K)$ in (\ref{Z-primes}) around this value, taking $n_{\max}$ to be of order $\ln{M}$, 
and then taking the limit $M \rightarrow \infty$ does lead to hyperuniformity
such that $Z(K) \sim K^2$ as $K\rightarrow 0$. 
\bigskip

\begin{figure}[bthp]
\includegraphics[  width=3in,keepaspectratio,clip=]{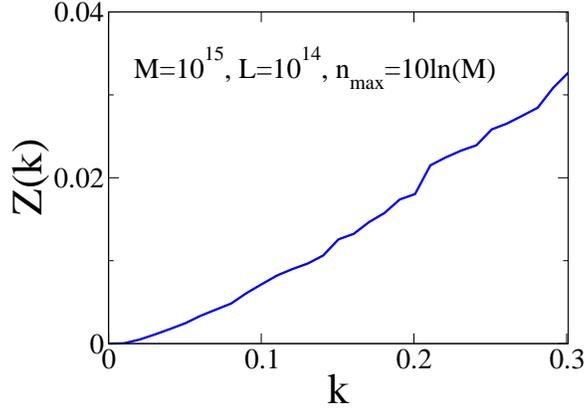}
\caption{The cumulative intensity function $Z(k)$ of the primes as obtained from (\ref{ZkPeak}). }
\label{Z-primes-plot}
\end{figure}

\noindent{\it Remark:} The result $Z(K) \sim K^2$ as $K\rightarrow 0$ implies that the primes fall within class II of
hyperuniform systems with a structure factor $S(k)$ that effectively is linear in $k$ as $k\rightarrow 0$ \cite{To18a}, so that the number variance $\sigma^2(R)$ scales logarithmically with $R$ in the large-$R$ limit.
Figure \ref{Z-primes-plot} depicts $Z(K)$ for the primes as determined from the first
line of relation (\ref{Z-primes}) in which $n_{\max}=10 \ln(M)$.

\section{$\tau$ Order Metric} \label{sec:tau}

We now study  the order metric $\tau$ in the discrete setting [c.f. (\ref{tau-discrete})] as a function of the system
size $L$ for the primes and the integer lattice.
\bigskip

\noindent{\sl Proposition 3:} The order metric $\tau$ of the 
primes as a function of the system size
$L$ in the considered interval $[M,M+L]$ (when divided by $\rho^2$) obeys the following scaling relation:
\begin{equation}
\tau/\rho^2 \sim c L,
\label{tau-L}
\end{equation} 
where $c$ is some constant.

This dependence of $\tau$ on $L$ can be understood in terms of the peaks in the structure factor via the definition (\ref{tau-discrete}). We have
\begin{equation}
\tau = \frac{1}{N_s} \sum_{j=1}^{N_s - 1} \left( S(j \pi/N_s) - 1+f\right)^2.
\end{equation}
We have seen how to estimate $S(\pi m/n)$ at a rational number in lowest terms: If the denominator is square-free, there will be a peak of height $N/\phi(n)^2$ while, if not, the structure factor will be very small. There could be a common factor between $j$ and $N_s$, so let
\begin{equation}
n = \frac{N_s}{\gcd(j,N_s)}.
\end{equation}
This makes $j/N_s = m/n$ in lowest terms with numerator $m = j/\gcd(j,N_s)$. Assume for simplicity that $N_s$ is square-free so that all of the resulting denominators $n$ will be square-free. Then we have
\begin{align*}
\tau &\sim \frac{1}{N_s} \sum_{n | N_s } \left( \frac{N}{\phi(n)^2} - 1 \right)^2 \#\{j \ \text{s.t.} \ N_s/\gcd(j,N_s) = n\} \\
&\sim \frac{N^2}{N_s} \sum_{n | N_s} \frac{1}{\phi(n)^4} \phi(n) \\
&\sim \rho^2 L \sum_{n | N_s} \frac{1}{\phi(n)^3}.
\end{align*}
The fact that the last sum is convergent proves the proposition.
Alternatively, one could use the expression (\ref{tau-discrete}) for $\tau$ and substitute the value predicted by Hardy-Littlewood for $g_2(2j)$.

In the case of the integer lattice, $\tau$ can be calculated directly from its definition (\ref{tau-discrete}). 
Since of the $N_s-1$ terms in the sum in (\ref{tau-discrete}), $1/f-1$ terms involve $S(k)$ at Bragg peak locations, while the remaining terms involved $S(k)$ values away from Bragg peaks, we have
\begin{equation}
\tau = \frac{1}{N_s} \left[\left(\frac{1}{f}-1\right) \left( N -1+f\right)^2 +  \left(N_s-\frac{1}{f}\right) \left( 0 - 1+f\right)^2\right].
\end{equation}
For sufficiently large systems, the contribution from the second term is negligible, and $N$ is much 
larger than $1-f$, yielding
\begin{equation}
\tau \approx \frac{1}{N_s} \left(\frac{1}{f}-1\right)  N^2 .
\end{equation}
Plugging in $N_s=L/2$, $f=2\rho$, and $N=\rho L$, one gets
\begin{align*}
\tau &\approx \frac{2}{L} \left(\frac{1}{2\rho}-1\right)  (\rho L)^2\\
&=L\rho(1-2\rho),
\end{align*}
which is consistent with Eq.~(\ref{tau-L}) for the primes when the density
is fixed. 

\begin{figure}[bthp]
\includegraphics[  width=3in]{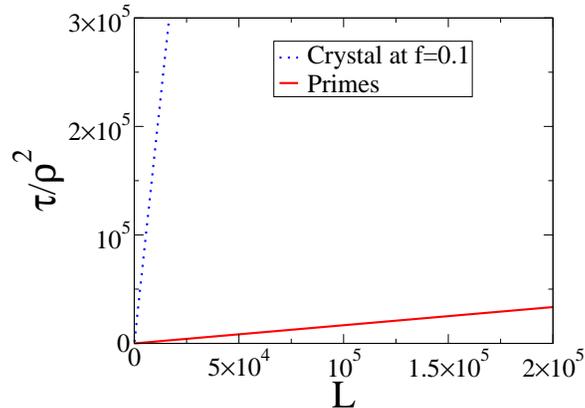}
\caption{ The order metric $\tau$, defined in Eq.~(\ref{tau-discrete}), in the discrete setting as a function of the system
size $L$ for the primes and integer lattice with filling fraction $f=0.1$ { obtained by direct simulations
from Eqs. (\ref{tau-discrete}) and (\ref{tau-L})}.}
\label{tau}
\end{figure}

We have numerically computed $\tau$ for the primes and integer
lattice by generating such configurations, sampling for their structure factors, and then 
computing their corresponding values of $\tau$ using relation (\ref{tau-discrete}).
For the primes, we take $M \approx 10^8$, and therefore $f =2\rho \approx 2/\ln(M)\approx 0.11$. For the integer lattice,
we take $f=0.1$. The constant $c$ is 18.00 for the integer
lattice and 0.1674 for the primes. These results are plotted in Fig. \ref{tau}.
Notice that
for a single configuration of an uncorrelated lattice-gas, $\tau$ does not grow with $L$ and converges in probability to a constant of order unity [see discussion under Eq. (\ref{tau-discrete})]. This means that the primes in these prescribed intervals are substantially  more ordered than the uncorrelated lattice gas
and appreciably less ordered than an integer lattice.

\begin{figure}[bthp]
\includegraphics[  width=4in]{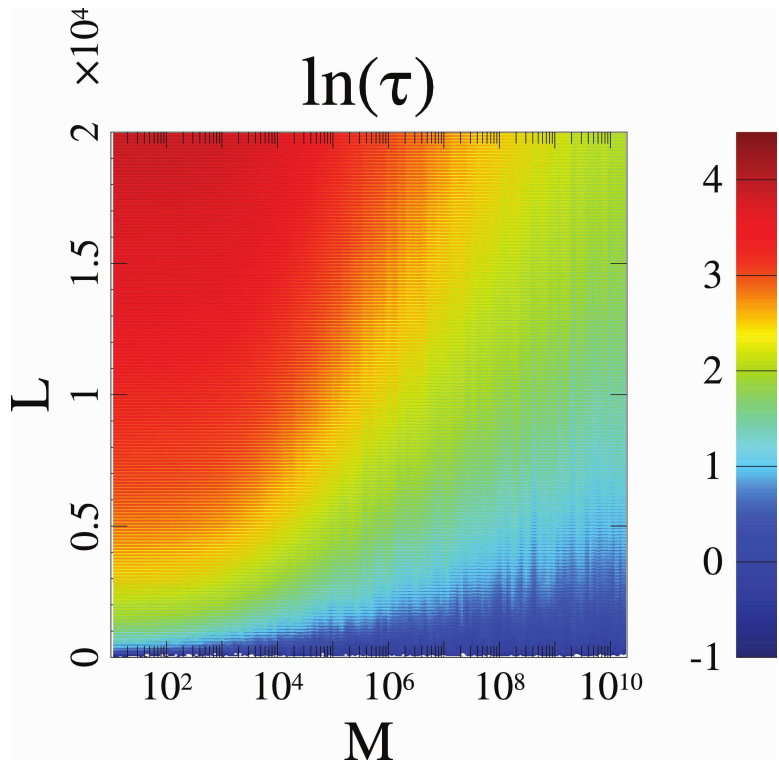}
\caption{Natural logarithm of the order metric $\tau$, defined in Eq.~(\ref{tau-discrete}), of prime numbers for $10<M<2\times 10^{10}$ and $8<L<2\times 10^4${ obtained from numerical simulations}.}
\label{tau-P}
\end{figure}

Now we study the $\tau$ order metric  of prime-number configurations with different
$M$ and $L$, obtained by computing their structure factors and evaluating (\ref{tau-discrete}).
As Fig.~\ref{tau-P} shows, the constant-$\tau$ level curves have the form $L \sim \ln^2(M)$, i.e., level curves appear
as quadratic curves in the log plot depicted  in Fig.~\ref{tau-P}.  
This follows from the fact that $\tau$ is proportional to $\rho^2 L$, the density $\rho$ being given by the prime number theorem as $1/\ln(M)$.
Thus, $L \sim \ln^2(M)$  is the boundary between regions where primes can be considered to be uncorrelated versus those where they exhibit correlations.

The regime $L \sim \ln(M)^2$ marks another important dividing line. Consider the Cram\'{e}r model where one replaces the number of primes from $M$ to $M+L$ by a sum of random variables taking the values 0 and 1 with probabilities chosen to match the prime number theorem. A sum of $L$ such random variables will have fluctuations on the order of $L^{1/2}$, by the Central Limit Theorem or by an elementary calculation, while the main term is $L/\ln(M)$. Thus, at least for the random model, $L \sim \ln(M)^2$ is a turning point between short intervals, where $L/\ln(M)^2 \rightarrow 0$ and the fluctuations overwhelm the prime number theorem, and longer intervals, where $L\ln(M)^2 \rightarrow \infty$ and the fluctuations can be neglected. Selberg \cite{selberg1943}, assuming the Riemann hypothesis, proved that for $L/\ln(M)^2 \rightarrow \infty$, the number of primes from $M$ to $M+L$ is $L/\ln(M)$ as predicted by the prime number theorem, except possibly for a sparse sequence of exceptional values of $M$. One might guess that this holds without exception, but Maier proved that there are infinitely many such $M$ \cite{maier1985}. For any power $c>1$, and setting $L=\ln(M)^c$, Maier proves
\begin{equation}
\limsup_{M \rightarrow \infty} \frac{N}{L/\ln(M)} > 1 > \liminf_{M \rightarrow \infty} \frac{N}{\ln(M)}.
\end{equation}
The behaviour of primes in intervals of this length is thus a significant departure from the random model.
For the uncorrelated regime in which Gallagher's results apply, $L\sim\ln(M)$, and $\tau$  will diminish as $M$ increases. 
As $L$ increases, prime-number configurations move from the uncorrelated regime ($\tau\sim 1$, $L \le \ln^2(M)$) to the {effective} limit-periodic
regime we studied in this paper ($\tau\sim L$, $L \propto M$), and then to the inhomogeneous regime where the density gradient is no longer
negligible (e.g., if $L \sim M^{1+\epsilon}$ and $\epsilon>0$.) In this last phase, the number
variance $\sigma^2(R)$ grows faster than the window volume (i.e., faster than $R$), which is the diametric opposite of hyperuniformity.

\section{Classification of the Primes} \label{sec:class}
\label{class}

Using Eq.~(\ref{S-p}), we have shown that the primes are like a limit-periodic point process, i.e.,
they are characterized by a structure factor $S(k)$ with dense Bragg peaks at certain
rational values of $k/\pi$ derived from an infinite union of periodic systems
in the limit as $M \to \infty$ with $L \sim \beta M$. This is similar to the structure factor of the limit-periodic 
systems; see Sec. \ref{doubling}. However, the primes show an erratic pattern of occupied
and unoccupied sites, very different from the predictable and orderly patterns of standard limit--periodic
systems.   Hence, the primes represent the first example of a point process that 
is {\it effectively} limit-periodic.

We have also demonstrated  that the primes are hyperuniform of class II.
Interestingly, this is precisely the same hyperuniformity class to which the (normalized) zeros
of the Riemann zeta function belong. 
However, the primes are substantially more ordered, having dense Bragg peaks instead of a continuous structure factor. As a result, one can only claim that $S(k) \rightarrow 0$ in a cumulative sense, unlike the case of the zeros.
While the Riemann zeros are disordered and hyperuniform,  the primes are effectively limit-periodic
and hence are characterized by order on all length scales. In terms of $\tau$, the primes are substantially  more ordered than the uncorrelated lattice gas 
and appreciably less ordered than an integer lattice, but similar in order to the period-doubling chain. It should not go unnoticed that
the lack of multiscale order in the Riemann zeros is reflected in the fact that $\tau$ is bounded in the large-$L$ limit. Indeed, assuming Montgomery's pair correlation conjecture, it converges to the constant 2/3
(which in the continuous setting should be compared to $\tau=0$ for the Poisson point process; see discussion
under Eq. (\ref{tau_c})]. This value 2/3 is what one finds by substituting the conjectured form for the structure factor 
(see Fig. \ref{nontrivial}) in the definition (\ref{tau_c}) of $\tau$ in the continuous setting 
and integrating $(k/(2\pi) - 1)^2$ over $0 < k < 2\pi$. In a system with multiscale order, such as the primes in dyadic intervals, $\tau$ diverges with $L$.

The condition that $L \sim M$, under which we have shown heretofore unanticipated order in the primes, is to be contrasted with the regime $L \sim \ln(M)$, in which the primes follow Gallagher's uncorrelated behavior. We have also shown that when $L$ grows faster than $M$, the primes enter the inhomogeneous regime where the density gradient is no longer negligible and hence where the limit-periodicity breaks down.

\section{Value Distribution of $S(k)$} \label{sec:values}

\begin{figure}[bthp]
\includegraphics[  width=3.8in,clip=]{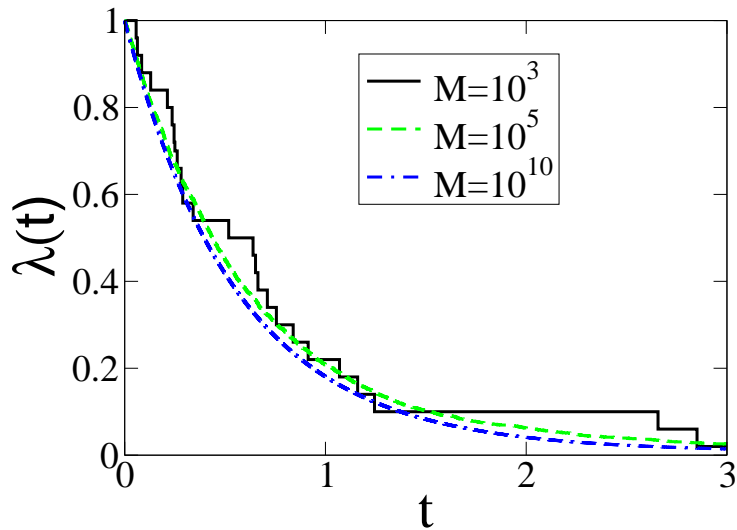}
\caption{The measure $\lambda(t)$ for relatively small range of values of the threshold $t$ for several
values of $M$. 
In all cases, $L=0.1 M$. The curves for different values of $M$ converges to an  
exponential curve well before the largest value of 
$M=10^{10}$ is attained. The latter case is well described by the exponential $\exp(-1.77989 t)$.}
\label{cdf-small}
\end{figure}

\begin{figure}[H]
\includegraphics[  width=4.2in,clip=]{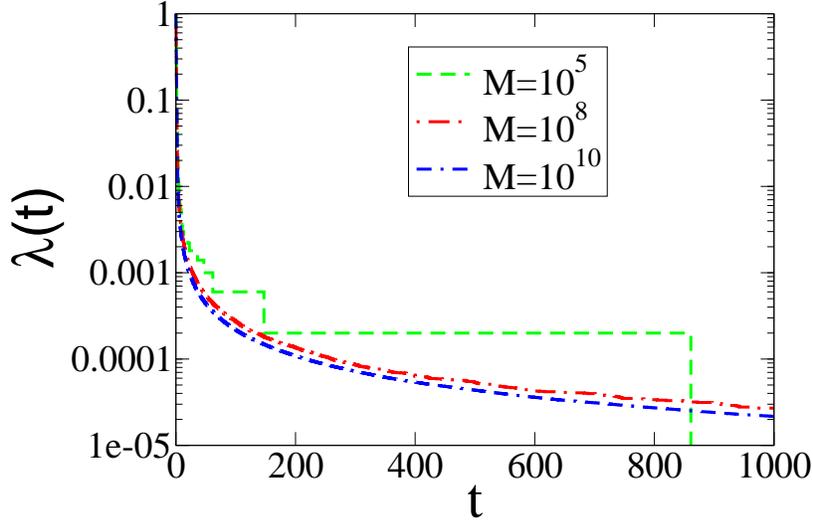}
\caption{The measure $\lambda(t)$ for relatively large range of values of the threshold $t$ for several
values of $M$. In all cases, $L=0.1 M$. Instead of the exponential curve found in the case shown in Fig. \ref{cdf-small}, here $\lambda(t)$
obeys an inverse power law. The curve for $M=10^{10}$ is well fit by the function $0.02174/t$, which we see satisfies the upper bound (\ref{eqn:cdf-bound}).}
\label{cdf-large}
\end{figure}

In this section, we study how frequently the structure factor $S(k)$ exceeds a given threshold $t$. Here, we write $S(k)$ for the structure factor of a finite system of $N$ primes $M < p < M+L$, namely
\begin{equation}
S(k) = \frac{1}{N} \left|\sum_p e^{ikp} \right|^2.
\end{equation}
Let 
\begin{equation} \label{eqn:cdf}
\lambda(t) = \frac{1}{\pi} | \{ 0 \leq k \leq \pi \ ; \ S(k) > t \} |
\end{equation}
be the measure of the set where $S(k) > t$, relative to the total length of the interval $0<k<\pi$. We think of $\lambda(t)$ as a cumulative distribution function (CDF). The quantity $\lambda(t)$ depends on $M$, not only $t$, but we suppress this dependence in the notation. There is also a related quantity which measures how often $S(k) > t$ while excluding forward scattering:
\begin{equation} \label{eqn:cdf-no-forward}
\lambda_{-}(t) = \left| \left\{ \frac{2\pi}{L} < k < \pi - \frac{2\pi}{L} \ ; \ S(k) > t \right\} \right|.
\end{equation}
We have
\begin{equation}
\lambda_{-}(t) \leq \lambda(t) \leq \lambda_{-}(t) + 2\pi/L,
\end{equation}
so the difference between $\lambda(t)$ and $\lambda_{-}(t)$ is negligible in the limit of large $L$. By orthogonality of the exponentials $e^{ikp}$, we have
\begin{equation} \label{eqn:parseval}
\int_0^{\pi} S(k) dk = \pi.
\end{equation}
Since $S(k) \geq 0$, it follows that
\begin{equation}
\pi = \int_0^{\pi} S(k) dk \geq \int_{ \{S(k)>t \} } S(k) dk \geq t \pi \lambda(t).
\end{equation}
Thus we have an upper bound
\begin{equation} \label{eqn:cdf-bound}
\lambda(t) \leq \frac{1}{t}.
\end{equation}
For $t>1$, this is an improvement over the trivial bound $\lambda(t) \leq 1$. For small $t$, Fig. \ref{cdf-small} suggests that $\lambda(t)$ is close to $e^{-c t}$ for some $c > 0$. If one imagines $S(k)$ as a small, noisy contribution on top of the peaks in its limiting form, this suggests that the noise follows an exponential distribution. On the other hand, Fig. \ref{cdf-large} suggests that the upper bound (\ref{eqn:cdf-bound}) is the correct order of magnitude of $\lambda(t)$ for large $t$. For $t$ on the order of $N$, the only way to have $S(k) > t$ is for $k$ to be close to a peak. This leads to the ``heavy-tailed'' power-law behavior illustrated in Fig. \ref{cdf-large}.

\section{Reconstruction of the Primes} \label{sec:recon}

It is noteworthy that we not only have an analytical formula  for the structure factor $S(k)$, which is related to the modulus of the complex 
density variable ${\tilde \n}(k)$, defined by (\ref{etak}), but also for the phase of ${\tilde \n}(k)$. The analytical expression
for ${\tilde \eta}(k)$ of the primes enables us to reconstruct, in principle,  a prime-number configuration within
an arbitrary interval $[M, M+L]$ by obtaining the inverse Fourier transform of ${\tilde \n}(k)$.  In this section, we will describe our
procedure to reconstruct prime-number configurations and report its accuracy.

We reconstruct a prime-number configuration in an interval $[M, M+L]$ by the following steps:
\begin{enumerate}
\item Calculate $N=(M+L)/\ln(M+L)-M/\ln(M)$. By the Prime Number Theorem, $N$ is the approximate number of primes in this interval.
\item Initialize ${\tilde \n}(k)$ at all $k\neq 0$ to be zero. Set ${\tilde \n}(0)=N$.
\item Find all odd square-free numbers $n \leq n_{\max}$.
\item For each $n$, find all integers $m$ such that $0< m<n$ and $m$ is co-prime with $n$.
\item For each $n$ and $m$, we need to reconstruct a peak at $k=m\pi/n$. For the purpose of reconstructing this peak, the prime-number configuration can be treated as a periodic system of period $2n$. We then find whether each number in 
the first period is prime, and calculate
\begin{equation}
C_1=\frac{N}{N_s}\frac{n}{\phi(n)} {\sum_{j}}^{\times} \exp(ik2j),
\end{equation}
which is the contribution to ${\tilde \eta}$ from the first $n$ sites. The index $j$ runs over numbers up to $N$ such that $\gcd(j,N)=1$.
\item If $L$ is divisible by $n$, then the peak at $k=m\pi/n$ should have infinitesimal width. We therefore increase ${\tilde \n}(m\pi/n)$ by $(L/n)C_1$, which  
is the predicted value of ${\tilde \rho}(m\pi/n)$.
\item If $L$ is not divisible by $n$, then the peak at $k=m\pi/n$ should have a finite width. We therefore have to increase ${\tilde \n}(k)$ of all $k$ points adjacent to $m\pi/n$ by
\begin{equation}
{\tilde \eta}(k)=C_1\left(\frac{1-f^{\left\lfloor L/2n \right\rfloor}}{1-f}\right),
\label{eq:Offpeak}
\end{equation} 
where $f=\exp(-i k 2n)$ is the phase change between the contribution to ${\tilde \eta}$ from the first $n$ sites and that from the second $n$ sites.
Here, the criteria for an ``adjacent'' $k$ point is that the absolute value of the result from Eq.~(\ref{eq:Offpeak}) is larger than $\sqrt{N}$.
\item Perform an inverse Fourier transform of ${\tilde \n}(k)$ to find $\n(r)$.
\end{enumerate}

If we had a completely accurate prediction of ${\tilde \n}(k)$, the resulting local density
$\n(r)$ would be exactly one for each prime number and exactly zero for each composite number. 
The predicted ${\tilde \n}(k)$ is not completely accurate. A major reason for the inaccuracy is that
there should be infinitely many peaks, but we only  
consider a finite number of peaks for which $n<n_{max}$.  Therefore, the resulting $\n(r)$ is not exactly zero or one. We find $N$ numbers with the highest predicted $\n(r)$ and predict those numbers to be prime. Thus another source of error is that we always predict $N$ primes, whereas the true number of primes may be more or less than the estimate from the Prime Number Theorem.

\begin{figure}[h]
\includegraphics[width=0.7\textwidth]{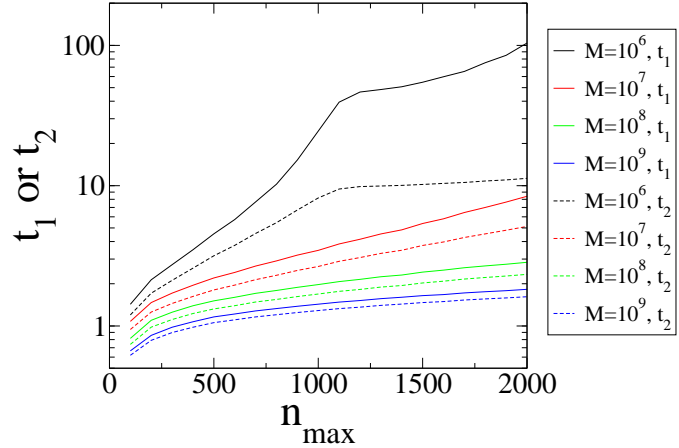}
\caption{Two measures of the accuracy of the predicted prime numbers, $t_1$ and $t_2$, defined by (\ref{1}) 
and (\ref{2}), respectively, versus $n_{max}$. }
\label{Prime_reconstruction}
\end{figure}

We have performed such reconstructions for $L=510510$ and several different $M$'s and $n_{max}$'s. 
Let $N_c$ be the number of correctly predicted primes, $N_i$ be the number of incorrect predictions, and 
$N_u$ be the number of un-predicted primes in this interval. We define the following ratios involving
these quantities to measure  the accuracy of the reconstruction procedure:
\begin{equation}
t_1=\frac{N_c}{N_i}, 
\label{1}
\end{equation}
\begin{equation}
t_2=\frac{N_c}{N_u}. 
\label{2}
\end{equation}
These two ratios versus $n_{max}$ are plotted in Fig.~\ref{Prime_reconstruction}. We see that for smaller $M$'s and larger $n_{max}$'s, this reconstruction procedure is highly accurate. When $M=10^{6}$ and $n_{max}=2000$, more than $99\%$ of the predicted prime numbers turn out to be correct. Unfortunately, as $M$ increases, the accuracy declines. For any $M$, increasing $n_{max}$ improves the accuracy, but with additional computational cost. We emphasize that an $n_{\max}$ of order $\ln{M}$ already yields nontrivial result. To find primes between $M$ and $2M$ naively by trial division, one would need to know the primes up to $\sqrt{M}$, much larger than our logarithmic $n_{\max}$. The fact that the first $\ln{M}$ primes are enough to predict primes in the long and distant interval from $M$ to $2M$, even imperfectly, is an interesting form of long-range order. Bertrand famously postulated that there is always a prime number between $M$ and $2M$, for any value of $M$. This was proved by Chebyshev using his explicit bounds towards the prime number theorem, but it remains a difficult algorithmic problem to find primes between $M$ and $2M$ quickly when $M$ is large. Our reconstruction procedure does give an algorithm
to do so. However, this is not a polynomial-time algorithm. The reason is that the Fast Fourier Transform scales linearithmically in the number of sites $L$, whereas the goal would be complexity polynomial in $\ln(L)$. Step (3) involves testing whether each number up to $n_{\max}$ is square-free, which up to logarithmic factors takes $O(n_{\max}^{3/2})$ steps (unless we had a faster way to test square-free than using trial division up to $\sqrt{n}$ to factor $n$). Step (4) takes $O(n_{\max}^2)$ steps, again up to logarithmic factors, using the Euclidean algorithm to find greatest common divisors quickly. Step (5) entails checking the primality of $n$ numbers, for $n$ up to $n_{\max}$ and computing a sum of length $N$ to find $C_1$. The primality testing can be done quickly and there is some overlap in the numbers that must be checked as $n$ varies. This step takes $O(n_{\max}^2 + N)$ operations. In step (7), we have to check whether each $k = m\pi/n$ with $n \leq n_{\max}$ is ``adjacent'', which involves $O(n_{\max}^2)$ cases. If $n_{\max}$ is negligible compared to $N$, then the slowest step is (8), where we perform the inverse Fourier transform of ${\tilde \n}(k)$. This step scales as $\mathcal{O}({L \ln{L}})$. Thus, at a computational cost of ${L \ln{L}}$ operations, we reconstruct the primes with imperfect accuracy. One could reconstruct with perfect accuracy by simply testing whether each of the $L$ numbers in the interval is prime. Each primality test takes at most on the order of $\ln(M)^6$ operations for AKS as adapted by Lenstra and Pomerance \cite{lenstra2002primality}, or $\ln(M)^4$ assuming the Riemann Hypothesis to run a deterministic Miller-Rabin test \cite{miller_1976}, \cite{rabin_1980}. Thus our method sacrifices accuracy in exchange for a running time faster by $\ln(L)^3$. 

It is interesting to note that when our reconstruction algorithm incorrectly predicts a composite number as prime, the composite number is usually ``almost prime'' in the sense that all of its prime factors are large. For example, our algorithm incorrectly predicted $1000733=809\times1237$ and $1001423=887\times1129$ as primes.


\section{Conclusions and Discussion}



The prime numbers display a range of behaviors depending on the interval under consideration. For dyadic intervals $[M,M+L]$ with $L$ comparable to $M$, we have found order 
across length scales, very different from the seeming randomness on display for smaller $L$. However, if $L$ were much larger than $M$, then one would reach the opposite 
conclusion of a non-hyperuniform system, purely because of the density gradient without reference to further properties specific to the primes.
The substantial order is reflected by the existence of dense Bragg peaks, a consequence of the {effective} limit-periodicity. The order metric $\tau$ gives a quantitative sense in which the
primes are substantially more ordered than the uncorrelated lattice gas but less ordered than an integer lattice.
When $L$ decreases, $\tau$ 
for the primes in this shorter interval becomes closer to that of an uncorrelated system with a transition visible at $L \sim \ln(M)^2$.
But the primes in dyadic intervals are hyperuniform, a fact which would seem almost unbelievable without the {effective} limit-periodic form for the
structure factor; see Eq. (\ref{S-p}).  Indeed,  the primes fall within the same broad hyperuniformity class 
as the Riemann zeta zeros (see Fig. \ref{nontrivial}), but are
substantially more ordered, having dense Bragg peaks instead of a continuous structure factor so that, unlike the Riemann zeros,
the order metric $\tau$ grows with $L$.  These peaks are located 
at rational wavenumbers $\pi m/n$ with odd, square-free denominator $n$, and were discovered numerically in \cite{Zh18}. They are explained in terms of the approximately equal distribution of prime numbers across residue classes modulo $2n$. Assuming the Hardy-Littlewood conjecture on prime pairs, the small ``diffuse part'' observed numerically in Ref. \cite{Zh18} is negligible in the infinite-system limit. The dense peaks exhibited by the primes are a feature shared with some recently studied quasicrystals, some also in class II  and others with even smaller density fluctuations, but these other examples have peaks  located at irrational wavenumbers \cite{Og17}. The primes are distinctive in being a superposition of periodic systems subject to irregularities in the distribution of occupied sites, which we call effective limit-periodicity.

The {effective} limit-periodic form of the structure factor allows one to predict the Hardy-Littlewood constants for the frequency with which $p$ and $p+r$ are both prime. 
Of course,  it remains an open problem to prove a lower bound establishing that there are infinitely many twin primes. As a more tractable open problem, we ask what further patterns can be found by considering three-particle and higher statistics, beyond what we could discern from the two-particle statistics $S(k)$ and $\tau$. This is related to the full Hardy-Littlewood $k$-tuples conjecture,  going beyond the case of $k=2$ considered here.

Since our analytical formula for the complex density variable ${\tilde \eta}(k)$,  defined by (\ref{etak}),
contains phase information, one can employ it to reconstruct a prime-number configuration within
an interval $[M, M+L]$ with $L \propto M$ by obtaining the inverse Fourier transform of ${\tilde \eta}(k)$.
This leads to an algorithm that enables the reconstruction of the primes in such intervals 
with high accuracy provided that $n_{\max}$ is sufficiently large and $M$ is not too large.



We are grateful to Peter Sarnak, Henry Cohn and Joshua Socolar for valuable discussions.
This work was supported in part by the National Science Foundation under Award No. DMR-1714722 and by the Natural Sciences and Engineering Research Council of Canada.

\appendix
\section{The Circle Method of Hardy-Littlewood}
For purposes of comparing our approach described in Sec. \ref{sec:circ} with the original approach of Hardy-Littlewood, 
we review the latter procedure. Their first objective is Goldbach's problem of writing a given $r$ as a \emph{sum} of primes $p_1 + \ldots + p_l$ rather than as a difference. Substantially the same analysis applies to the equations $p_1 + p_2 = n$ and $p_2 - p_1 = r$. These form an example of what Hardy-Littlewood call \emph{conjugate} problems. Instead of a trigonometric polynomial, their generating function is
\begin{equation}
f(x) = \sum_p \ln(p) x^p
\end{equation}
where the sum is over all primes, not just those in a finite interval. For the sum to converge, we must have $|x|<1$. It is $f(x)^2$ that features in Goldbach's problem, versus $|f(x)|^2$ for prime pairs.
The main theorem of \cite{Ha23} is that, assuming a strong hypothesis on zeros of Dirichlet $L$-functions, every sufficiently large odd number is a sum of three primes. One can fairly take this hypothesis to be that all zeros of all Dirichlet $L$-functions lie on the critical line, namely the generalized Riemann hypothesis, but Hardy-Littlewood in fact work under the weaker hypothesis of a certain zero-free strip. The role of the generalized Riemann hypothesis is to ensure that primes are evenly distributed in arithmetic progressions to such large moduli as $q = \sqrt{n}$, which occur in the Farey dissection used by Hardy-Littlewood. In the unconditional work of Vinogradov, the denominators grow only as powers of $\ln(n)$ instead of $\sqrt{n}$. With or without the Riemann hypothesis, the method fails for $l=2$, and this is why we do not have a proof of (\ref{S-p}). The main term obtained by Hardy-Littlewood amounts to $n^{l-1}$, with errors no greater than $n^{l/2 + 1/4 + \epsilon}$ in order of magnitude.
Thus the estimate succeeds if $l-1 > l/2 + 1/4$, and so for all larger $l \geq 3$ but not $l=2$.

The circle method, aptly named, begins with an integral over a circle $|x| = R$ in the complex plane:
\[
\frac{1}{2\pi } \int_0^{2\pi} f(Re^{-i\psi}) f(Re^{-i\psi}) e^{i r \psi} d\psi = R^r \sum_p a_p R^{2p}
\]
where $a_p = \ln(p)\ln(p+r)$ if both $p$ and $p_r$ are prime, and $a_p = 0$ otherwise. As $R \rightarrow 1$, this gives a count of the number of prime pairs. Note that $f(Re^{i\psi})f(Re^{-i\psi}) = |f(x)|^2$. To estimate the integral, Hardy and Littlewood divide the circle into Farey arcs. The arc $\xi_{a,q}$ consists of those $x = Re^{i\psi}$ whose argument lies in a particular interval
\[
\frac{2\pi a}{q} - \theta_{a,q}' < \psi < \frac{2\pi a}{q} + \theta_{a,q}.
\]
The values of $\theta_{a,q}$ and $\theta_{a,q}'$ are such that these arcs cover the circle without overlap, and the length of $\xi_{a,q}$ is of order $(qN)^{-1}$, $N$ being the last stage to which one carries out the Farey dissection. This $N$ is what we call $q_{\max}$. Writing the radius $R$ in the form $R=e^{-1/n}$, Hardy and Littlewood choose $N = \lfloor \sqrt{n} \rfloor$. On the arc $\xi_{a,q}$, $f(x)$ is approximated by
\begin{equation} \label{eq:lemma9}
f(x) \approx \frac{\mu(q)}{\phi(q)} \frac{1}{1/n - i (\psi - 2\pi a/q )}
\end{equation}
with an estimate on the error given by their Lemma 9 p. 19 \cite{Ha23}. The error is essentially the square root of the main term, provided $q \leq \sqrt{n}$ and assuming the Riemann hypothesis. Changing variables to $u = \psi - 2\pi a/q$ leads to the simple expression
\[
|f(x)|^2 \approx \frac{\mu(q)^2}{\phi(q)^2} \frac{1}{1/n^2 + u^2}.
\]
Integrating over $\xi_{a,q}$ then gives
\[
\int_{\xi_{a,q}} f(Re^{i\psi})f(Re^{-i\psi}) e^{ir\psi} d\psi \approx \frac{\mu(q)^2}{\phi(q)^2} e^{2\pi iar/q} \int_{-\theta_{a,q}'}^{\theta_{a,q}} \frac{e^{iru}}{u^2 +1/n^2} du
\]
For large $n$, we may allow the further approximation that
\[
\int_{\theta'}^{\theta} \frac{e^{iru}}{u^2 + 1/n^2} du = n \int_{-n\theta'}^{n\theta} \frac{e^{iv (r/n)}}{v^2+1}dv \sim \pi n.
\]
The factor $\pi$ is obtained by replacing the integral with $\int_{-\infty}^{\infty} (v^2+1)^{-1} dv$. Note that $\theta$ and $\theta'$ are of order $(qN)^{-1} \sim q^{-1} n^{-1/2}$, so that $n\theta \asymp q^{-1} \sqrt{n} \rightarrow \infty$, while the exponential $e^{i v r/n}$ is close to 1 for $n$ large, 

Combining all of the arcs $\xi_{a,q}$ such that $q \leq \sqrt{n}$ and $a \leq q$, with no common factor between $a$ and $q$, we find that
\[
R^r \sum_p a_p R^{2p} \sim \frac{1}{2\pi} \sum_q \sum_a \frac{\mu(q)^2}{\phi(q)^2} e^{2\pi i ar/q} \pi n = \frac{n}{2} \mathfrak{S}_2
\]
To write this in terms of $R = e^{-1/n}$, note that $1 - R^2 = 1 - e^{-2/n} \sim 2/n$.
Therefore
\[
\sum_p a_p R^{2p} \sim \frac{1}{1-R^2} \mathfrak{S}_2.
\]
Hardy and Littlewood employ a Tauberian argument to deduce from this behaviour as $R \rightarrow 1$ that
\[
\sum_{p<n} a_p \sim n \mathfrak{S}_2.
\]
The ability to skip this step is one of the technical advantages of working with a finite generating function such as $S(k)$ instead of $|f(x)|^2$.
Finally, removing the logarithmic weights in $a_p$, Hardy and Littlewood formulate their Conjecture B: As $n \rightarrow \infty$, the number of prime pairs $p, p+r$ less than $n$ is asymptotic to
\[
\frac{n}{\ln(n)^2} 2 C_2  \prod_{p | r} \left( \frac{p-1}{p-2} \right)
\]
where $C_2$ is their twin primes constant:
\[
C_2 = \prod_{p > 2} (1 - 1/(p-1)^2 ).
\]
Hardy and Littlewood rewrite $\mathfrak{S}_2$ as a product over primes using their Lemma 12 ``Summation of the singular series" (p. 27 of \cite{Ha23}):
\[
\sum_{q=1}^{\infty} \left( \frac{\mu(q)}{\phi(q)} \right)^l c_q (-r) = 2 \prod_{\varpi = 3}^{\infty} \left( 1 - \frac{(-1)^l}{(\varpi - 1)^l} \right) \prod_{p | r} \left( \frac{(p-1)^l +(-1)^l (p-1)}{(p-1)^l - (-1)^l} \right)
\]
assuming $r$ and $l$ are of the same parity, whereas the sum is 0 if $r$ and $l$ have opposite parity. In this notation, $p$ refers to an odd prime divisor of $r$ while $\varpi$ is any odd prime, and $c_q(-r)$ is Ramanujan's sum $\sum_a e^{-2\pi i ar/q}$ over $a$ modulo $q$ with no factor in common with $q$. When $l=2$, the product over primes becomes
\[
2 \prod_{\varpi = 3}^{\infty} \left( 1 - \frac{1}{(\varpi - 1)^2} \right) \prod_{p | r} \left( \frac{(p-1)^2 + (p-1)}{(p-1)^2- 1} \right) = 2 C_2 \prod_{p | r} \left( \frac{p-1}{p-2} \right)
\]
hence the value in Conjecture B.
To see why this infinite product agrees with the sum, imagine multiplying together several terms from the product. Note that all primes $p$ and $\varpi$ in this formula are odd, the factor 2 having already accounted for $p=2$, which divides $r$. 
Each $p$ dividing $r$ introduces a factor $(1 + 1/(p-1))$, since we must incorporate the factor in $C_2$ from $\varpi = p$. For $p$ not dividing $r$, the factor is $(1-1/(p-1)^2)$. 
The result is a sum over all $n = p_1 \cdots p_t$, namely odd square-free numbers. %




\begin{thebibliography}{10}

\bibitem{atkin1993elliptic}
A.~O.~L. Atkin and F.~Morain, \emph{Elliptic curves and primality proving},
  Math. Computation \textbf{61} (1993), no.~203, 29--68.

\bibitem{Ba11}
M.~Baake and U.~Grimm, \emph{Diffraction of limit periodic point sets}, Phil.
  Mag. \textbf{91} (2011), 2661--2670.

\bibitem{baillie1980lucas}
R.~Baillie and S.~S. Wagstaff, \emph{Lucas pseudoprimes}, Math. Computation
  \textbf{35} (1980), no.~152, 1391--1417.

\bibitem{Ba08}
R.~D. Batten, F.~H. Stillinger, and S.~Torquato, \emph{Classical disordered
  ground states: {S}uper-ideal gases, and stealth and equi-luminous materials},
  J. Appl. Phys. \textbf{104} (2008), 033504.

\bibitem{Bom86}
E.~Bombieri and J.~E. Taylor, \emph{Which distributions of matter diffract?
  {A}n initial investigation}, J. de Phys. Coll. \textbf{47} (1986), C3--19.

\bibitem{cherwell}
Lord Cherwell, \emph{On the distribution of the intervals between prime
  numbers}, Quart. J. Math. \textbf{17} (1946), 46--62.

\bibitem{dahmen_2001}
S.~R. Dahmen, S.~D. Prado, and T.~Stuermer-Daitx, \emph{Similarity in the
  statistics of prime numbers}, Physica A \textbf{296} (2001), 523--528.

\bibitem{Da80}
H.~Davenport, \emph{Multiplicative number theory}, vol.~74, Springer, New York,
  1980.

\bibitem{Di18}
R.~A. DiStasio, G.~Zhang, F.~H. Stillinger, and S.~Torquato, \emph{Rational
  design of stealthy hyperuniform patterns with tunable order}, Phys. Rev. E
  \textbf{97} (2018), 023311.

\bibitem{Do05d}
A.~Donev, F.~H. Stillinger, and S.~Torquato, \emph{Unexpected density
  fluctuations in disordered jammed hard-sphere packings}, Phys. Rev. Lett.
  \textbf{95} (2005), 090604.

\bibitem{Dy62a}
F.~J. Dyson, \emph{Statistical theory of the energy levels of complex systems.
  {I}}, J. Math. Phys. \textbf{3} (1962), 140--156.

\bibitem{Fl09b}
M.~Florescu, S.~Torquato, and P.~J. Steinhardt, \emph{Designer disordered
  materials with large complete photonic band gaps}, Proc. Nat. Acad. Sci.
  \textbf{106} (2009), 20658--20663.

\bibitem{gallagher1976distribution}
P.~X. Gallagher, \emph{On the distribution of primes in short intervals},
  Mathematika \textbf{23} (1976), 4--9.

\bibitem{Go17}
T.~Goldfriend, H.~Diamant, and T.~A. Witten, \emph{Screening, hyperuniformity,
  and instability in the sedimentation of irregular objects}, Phys. Rev. Lett.
  \textbf{118} (2017), 158005.

\bibitem{Gr95}
A.~Granville, \emph{Harald {C}ram{\'e}r and the distribution of prime numbers},
  Scand. Actuarial J. \textbf{1995} (1995), 12--28.

\bibitem{green_2008_2}
B.~Green and T.~Tao, \emph{The mobius function is asymptotically orthogonal to
  nilsequences}, arXiv:0807.1736 [math.NT] (2008).

\bibitem{green_2008}
\bysame, \emph{The primes contain arbitrarily long arithmetic progressions},
  Annals of Math. \textbf{167} (2008).

\bibitem{green_2006}
B.~Green and T.~Tao, \emph{Linear equations in primes}, Annals Math. (2010),
  1753--1850.

\bibitem{hadamard1896distribution}
J.~Hadamard, \emph{Sur la distribution des z{\'e}ros de la fonction $\zeta(s)$
  et ses cons{\'e}quences arithm{\'e}tiques}, Bulletin de la Societ{\'e}
  mathematique de France \textbf{24} (1896), 199--220.

\bibitem{Ha23}
G.~H. Hardy and J.~E. Littlewood, \emph{Some problems of ‘partitio
  numerorum’; {III}: On the expression of a number as a sum of primes}, Acta
  Mathematica \textbf{44} (1923), 1--70.

\bibitem{He17b}
D.~Hexner, P.~M. Chaikin, and D.~Levine, \emph{Enhanced hyperuniformity from
  random reorganization}, Proc. Nat. Acad. Sci. \textbf{114} (2017),
  4294–--4299.

\bibitem{Iw04}
H.~Iwaniec and E.~Kowalski, \emph{Analytic number theory}, vol.~53, American
  Mathematical Society, 2004.

\bibitem{Ja15}
R.~L. Jack, I.~R. Thompson, and P.~Sollich, \emph{Hyperuniformity and phase
  separation in biased ensembles of trajectories for diffusive systems}, Phys.
  Rev. Lett. \textbf{114} (2015), 060601.

\bibitem{Ji14}
Y.~Jiao, T.~Lau, H.~Hatzikirou, M.~Meyer-Hermann, J.~C. Corbo, and S.~Torquato,
  \emph{Avian photoreceptor patterns represent a disordered hyperuniform
  solution to a multiscale packing problem}, Phys. Rev. E \textbf{89} (2014),
  022721.

\bibitem{lenstra2002primality}
Hendrik~W Lenstra~Jr and Carl Pomerance, \emph{Primality testing with gaussian
  periods}, Foundations of Software Technology and Theoretical Computer
  Science, 2002, p.~1.

\bibitem{Lev86}
D.~Levine and P.~J Steinhardt, \emph{Quasicrystals. {I}. {D}efinition and
  structure}, Phys. Rev. B \textbf{34} (1986), 596.

\bibitem{Ma16}
T.~Ma, H.~Guerboukha, M.~Girard, A.~D. Squires, R.~A. Lewis, and
  M.~Skorobogatiy, \emph{3d printed hollow-core terahertz optical waveguides
  with hyperuniform disordered dielectric reflectors}, Adv. Optical Mater.
  \textbf{4} (2016), 2085--2094.

\bibitem{maier1985}
H.~Maier, \emph{Primes in short intervals}, Michigan Math. J. \textbf{32}
  (1985), 221--225.

\bibitem{martelli_2013}
F.~Martelli, \emph{Dealing with primes i.: On the goldbach conjecture},
  arXiv:1309.5895 [math.NT] (2013).

\bibitem{Ma15}
A.~Mayer, V.~Balasubramanian, T.~Mora, and A.~M. Walczak, \emph{How a
  well-adapted immune system is organized}, Proc. Nat. Acad. Sci. \textbf{112}
  (2015), 5950--5955.

\bibitem{maynard_2015}
J.~Maynard, \emph{Small gaps between primes}, Annals of Math. \textbf{181}
  (2015), 383--413.

\bibitem{Me91}
M.~L. Mehta, \emph{Random matrices}, Academic Press, New York, 1991.

\bibitem{miller_1976}
G.~Miller, \emph{Riemann's hypothesis and tests for primality}, Journal of
  Computer and System Sciences \textbf{13} (1976), 300--317.

\bibitem{Mon73}
H.~L. Montgomery, \emph{The pair correlation of zeros of the zeta function},
  Amer. Math. Soc. (1973), 181--193.

\bibitem{Mo04}
H.~L. Montgomery and K.~Soundararajan, \emph{Primes in short intervals}, Comm.
  Math. Phys. \textbf{252} (2004), 589--617.

\bibitem{lemke_2016}
R.~J.~Lemke Oliver and K.~Soundararajan, \emph{Unexpected biases in the
  distribution of consecutive primes}, arXiv:1603.03720 [math.NT] (2016).

\bibitem{Og17}
E.~C. {O{\u g}uz}, J.~E.~S. {Socolar}, P.~J. {Steinhardt}, and S.~{Torquato},
  \emph{{Hyperuniformity of Quasicrystals}}, Phys. Rev. B \textbf{95} (2017),
  054119.

\bibitem{pintz2007}
J.~Pintz, \emph{{C}ram\'{e}r vs. {C}ram\'{e}r. {O}n {C}ram\'{e}r's
  probabilistic model for primes}, Functiones et Approximatio (2007), 361--376.

\bibitem{pomerance1980pseudoprimes}
C.~Pomerance, J.~L. Selfridge, and S.~S. Wagstaff, \emph{The pseudoprimes to
  25· 10?}, Math. Computation \textbf{35} (1980), no.~151, 1003--1026.

\bibitem{rabin_1980}
M.~Rabin, \emph{Probabilistic algorithm for testing primality}, Journal of
  Number Theory \textbf{12} (1980), 128--138.

\bibitem{rubinstein_1994}
M.~Rubinstein and P.~Sarnak, \emph{Chebyshev’s bias}, Experiment. Math.
  \textbf{3} (1994), 173--197.

\bibitem{Rud96}
Z.~Rudnick and P.~Sarnak, \emph{Zeros of principal {$L$}-functions and random
  matrix theory}, Duke Math. J. \textbf{81} (1996), 269--322.

\bibitem{selberg1943}
A.~Selberg, \emph{On the normal density of primes in small intervals, and the
  difference between consecutive primes}, Archiv for Mathematik og
  Naturvidenskab \textbf{B. 47, No. 6} (1943), 87--105.

\bibitem{tao_2011}
T.~Tao and T.~Ziegler, \emph{The inverse conjecture for the gowers norm over
  finite fields in low characteristic}, Annals of Combinatorics \textbf{16}
  (2012), no.~1, 121--188.

\bibitem{Te95}
G.~Tenenbaum, \emph{Introduction to analytic and probabilistic number theory},
  vol.~46, Cambridge University Press, Cambridge, 1995.

\bibitem{To18a}
S.~Torquato, \emph{Hyperuniform states of matter}, Physics Reports,
  arXiv:1801.06924 (2018).

\bibitem{To08c}
S.~Torquato, A.~Scardicchio, and C.~E. Zachary, \emph{Point processes in
  arbitrary dimension from {F}ermionic gases, random matrix theory, and number
  theory}, J. Stat. Mech.: Theory Exp. (2008), P11019.

\bibitem{To03a}
S.~Torquato and F.~H. Stillinger, \emph{Local density fluctuations,
  hyperuniform systems, and order metrics}, Phys. Rev. E \textbf{68} (2003),
  041113.

\bibitem{To18b}
S.~Torquato, G.~Zhang, and M.~de~Courcy-Ireland, \emph{Uncovering multiscale
  order in the prime numbers via scattering},  (2018).

\bibitem{To15}
S.~{Torquato}, G.~{Zhang}, and F.~H. {Stillinger}, \emph{{Ensemble Theory for
  Stealthy Hyperuniform Disordered Ground States}}, Phys. Rev. X \textbf{5}
  (2015), 021020.

\bibitem{vaughan97}
R.~C. Vaughan, \emph{The {H}ardy-{L}ittlewood method}, vol. 125, Cambridge
  University Press, Cambridge, 1997.

\bibitem{vinogradov_1937}
I.~M. Vinogradov, \emph{The method of trigonometrical sums in the theory of
  numbers (russian)}, Trav. Inst. Math. Steklo \textbf{10} (1937).

\bibitem{wolf_1996}
M.~Wolflf, \emph{Unexpected regularities in the distribution of prime numbers},
  Proc. of the 8th Joint EPS - APS Int. Conf. Physics Computing (1996).

\bibitem{Za11a}
C.~E. Zachary, Y.~Jiao, and S.~Torquato, \emph{Hyperuniform long-range
  correlations are a signature of disordered jammed hard-particle packings},
  Phys. Rev. Lett. \textbf{106} (2011), 178001.

\bibitem{Za09}
C.~E. Zachary and S.~Torquato, \emph{Hyperuniformity in point patterns and
  two-phase heterogeneous media}, J. Stat. Mech.: Theory \& Exp. (2009),
  P12015.

\bibitem{Zh18}
G.~Zhang, F.~Martelli, and S.~Torquato, \emph{Structure factor of the primes},
  J. Phys. A: Math. \& Gen. \textbf{51} (2018), 115001.

\bibitem{zhang_2014}
Y.~Zhang, \emph{Bounded gaps between primes}, Annals of Math. \textbf{3}
  (2014), 1121--1174.

\end{thebibliography}


\providecommand{\bysame}{\leavevmode\hbox to3em{\hrulefill}\thinspace}
\providecommand{\MR}{\relax\ifhmode\unskip\space\fi MR }
\providecommand{\MRhref}[2]{%
  \href{http://www.ams.org/mathscinet-getitem?mr=#1}{#2}
}
\providecommand{\href}[2]{#2}


\end{document}